\documentclass[12pt]{article}

\newtheorem{thm}{Theorem}[section]
\newtheorem{prop}[thm]{Proposition}
\newtheorem{lem}[thm]{Lemma}
\newtheorem{dfn}[thm]{Definition}

\newtheorem{con}[thm]{Conjecture}

\newcommand{\qed}{\hfill \fbox{}\medskip}
\newcommand{\proof}{\medskip\noindent{\bf Proof.}\quad }

\title{A Commutative Family of Integral Transformations
and Basic Hypergeometric Series. II. 
Eigenfunctions and Quasi-Eigenfunctions}

\author{Jun'ichi Shiraishi\\
\\
{\it Graduate School of 
Mathematical Science, }\\ {\it University of Tokyo, Tokyo, Japan}}

\date{}

\begin{document}

\maketitle

\maketitle

\begin{abstract}
A series of conjectures is obtained as further investigation of the
integral transformation $I(\alpha)$ introduced in the previous paper.
A Macdonald-type difference operator $D$ is introduced. 
It is conjectured that $D$ and $I(\alpha)$ are commutative with each other. 
Studying the series for the eigenfunctions under termination conditions,
it is observed that a deformed Weyl group action appears as
a hidden symmetry.
An infinite product formula for the eigenfunction is 
found for a spacial case of parameters.
A one parameter family of
hypergeometric-type series $F(\alpha)$ is introduced.
The series $F(\alpha)$ is caracterized by 
a covariant transformation property 
$I(\alpha q^{-1} t)\cdot F(\alpha)=F(\alpha q^{-1} t)$
and a certain initial condition given at $\alpha=t^{1/2}$.
We call $F(\alpha)$ the `quasi-eigenfunction' for short.
A class of infinite product-type expressions are conjectured
for $F(\alpha)$ at the special points 
$\alpha=-t^{1/2}$, $\alpha=t$, $\alpha=\pm q^{1/2}t^{1/2}$,
and  $\alpha=\pm q^{\ell}t^{1/2}$ ($\ell=1,2,3,\cdots$).
\end{abstract}

\section{Introduction}
This is the second paper of a series in which 
a commutative family generated by an integral transformation is studied.
In the last paper \cite{S1}, we have introduced the integral transformation 
$I(\alpha)$ and conjectured the commutativity $[I(\alpha),I(\beta)]=0$.
Our arguments there were based on 
the basic properties for the the eigenfunctions and
the explicit formulas.
The main purpose of the present paper is
to state a series of conjectures which 
gives us several useful
descriptions for the structure of the operator $I(\alpha)$
and its eigenfunctions. We consider these as 
preparations for our application of $I(\alpha)$
to the investigation of the vertex operators for the eiget-vertex model 
given in the next paper \cite{S3}.
As for the eight-vertex model and the corner transfer matrix approach 
to the correlation functions, we refer the readers to 
\cite{B1,B2,B} and \cite{FIJKMY1,FIJKMY2}.
\bigskip

The integral transformation 
$I(\alpha)$ acting on the space of formal power series 
${\cal F}_n={\bf C}[[\zeta_2/\zeta_1,\cdots,\zeta_n/\zeta_{n-1}]]$
was introduced in the first paper \cite{S1} (see Definition 1.1 of \cite{S1})
as follows:
\begin{eqnarray}
&&I(\alpha)  f(\zeta_1,\cdots,\zeta_n) \nonumber\\
&=&
\prod_{i=1}^n \left( 
{ (qt^{-1} ;q)_\infty \over 
(\alpha s_i^{-1} qt^{-1} ;q)_\infty}
{ (q;q)_\infty \over 
(\alpha^{-1}s_i q;q)_\infty}
\right)\\
&&
\times\prod_{i<j} h(\zeta_j/\zeta_i)
\oint_{C_1}\cdots\oint_{C_n} {d\xi_1 \over 2\pi i\xi_1}\cdots
{d\xi_n \over 2\pi i\xi_n}
\prod_{i=1}^n
{\Theta_{q}(\alpha s_i^{-1} q^{1\over 2}t^{-{1\over 2}} \zeta_i/\xi_i)\over 
\Theta_{q}(q^{1\over 2}t^{-{1\over 2}}\zeta_i/\xi_i)} \nonumber\\
&&
\times
\prod_{k=1}^n
\left[
\prod_{i<k} g(\zeta_k/\xi_i)
\prod_{j\geq k} g(\xi_j/\zeta_k)
\right]
f(\xi_1,\cdots,\xi_n), \nonumber
\end{eqnarray}
where the integration contours $C_i$ are given by the conditions $|\zeta_i/\xi_i|=1$,
and the functions $h(\zeta)$ and $g(\zeta)$ are given by
\begin{eqnarray}
h(\zeta)=(1-\zeta)
{(q t^{-1} \zeta;q)_\infty \over
(t\zeta;q)_\infty },\qquad
g(\zeta)=
{(q^{1\over 2}t^{1\over 2} \zeta;q)_\infty \over
(q^{1\over 2}t^{-{1\over 2}} \zeta;q)_\infty }.
\end{eqnarray}
The main statement there was 
the conjecture for the commutativity (see Conjecture 1.2 in \cite{S1}):
\begin{eqnarray}
I(\alpha) I(\beta)=I(\beta) I(\alpha). \label{II=II}
\end{eqnarray}
A proof of this for the case $n=2$ was obtained by using some 
summation and transformation formulas for the basic hypergeometric series.
In this paper, we continue our study on the 
operator $I(\alpha)$ and try to obtain 
better understanding. 
\medskip

This paper consist of several conjectures 
(proved for some cases) and definitions summarized as follows.
\begin{enumerate}
\item Commutativity between the integral operator $I(\alpha)$ and a Macdonald-type difference operator $D$ acting on ${\cal F}_n$.
\item Termination of the eigenfunctions of $I(\alpha)$ and a hidden 
symmetry of a Weyl group action.
\item Infinite product formulas for the eigenfunctions of $I(\alpha)$
for a particular case of parameters.
\item Introduction of the `quasi-eigenfunction' $F(\alpha)$ which
is defined by a certain transformation property 
with respect to the action of $I(\alpha)$ and a certain initial condition.
\item
A class of infinite product(-type) formulas for $F(\alpha)$
for particular values of $\alpha$.
\end{enumerate}

The Jack and the Macdonald symmetric polynomials \cite{Mac}
and the raising operators were studied in \cite{AMOS}\cite{AOS}. 
The raising operator can be characterized as
an eigenfunction of a Macdonald-type difference operator $D$.
As for the detail,
we refer the readers to Appendix A 
and the lecture note \cite{S-Lec}. 
Note that the solution can be
regarded as a basic analogue of 
the Heckman-Opdam hypergeometric function \cite{HO1}.
While obtaining the explicit formulas for the eigenfunctions of $I(\alpha)$,
the author realized that $I(\alpha)$ and $D$ share exactly the same 
eigenfunctions, at least for small $n$. This coincidence seems quite mysterious
and meaningful, because
it suggests that there exists a profound relationship between
the eight-vertex model and the Macdonald polynomials.
\medskip

In \cite{S1}, explicit formulas for $I(\alpha)$ for small $n$ were studied.
(See Theorem 2.1, Conjecture 3.3 and Proposition 3.4 in the first paper \cite{S1}.)
We observed that the hypergeometric-type series expressions
obtained there 
gives us an efficient way to organize the series for the eigenfunctions.
These explicit expressions allow us to find 
the termination conditions for the series. 
One may observe that the resulting polynomials 
have a good symmetry.
We will introduce a certain Weyl group action on the space of
polynomials, 
to explain this hidden symmetry for the eigenfunctions.
\medskip

We observe that the eigenfunction of $I(\alpha)$ 
can be written as an infinite product, if we specialize 
the parameter $s_i$'s as
$(s_1,s_2\cdots,s_n)=(1,t,\cdots,t^{n-1})$.
This observation comes from the explicit formulas 
for small $n$. If we assume the above stated 
commutativity $[I(\alpha),D]=0$,
this factorization can be shown for general $n$.
\medskip

Next, we introduce another class of 
`multi variable hypergeometric series' depending on a parameter $\alpha$.
Setting the parameters as
$(s_1,s_2\cdots,s_n)=(1,1,\cdots,1)$,
we characterize a continuous family of series $F(\alpha)$ by
imposing the following two conditions:
\begin{eqnarray}
{\rm (I)}&&
I(\alpha q^{-1}t ) \cdot F(\alpha)=F(\alpha q^{-1}t),\\
{\rm (II)}&&
F(t^{1/2})
=\prod_{1\leq i<j\leq n}(1-\zeta_j/\zeta_i)
{(q t^{-1/2}\zeta_j/\zeta_i;q)_\infty \over 
(t^{1/2}\zeta_j/\zeta_i;q)_\infty }.
\end{eqnarray}
We call this function $F(\alpha)$ `quasi-eigenfunction' for short.
In the next paper \cite{S3}, it will be explained that 
this definition for $F(\alpha)$ naturally emerges 
in the context of the study of the eight-vertex mode 
based on the corner transfer matrix method.
\medskip

It will be argued that 
a class of infinite product(-type)
expressions can be found for the quasi-eigenfunctions $F(\alpha)$
at particular values for $\alpha$.
Note that some of these product formulas 
were observed in the previous paper \cite{S}, 
while we studied the eight-vertex model for 
the parameters $p^{1/2}_{\rm 8v}=q^{3/2}_{\rm 8v}, -q^2_{\rm 8v}$ and $q^3_{\rm 8v}$. 
(Note that the basic parameters $p_{\rm 8v}$ and $q_{\rm 8v}$ used in \cite{S} are
related with the ones in this paper as
$p^{1/2}_{\rm 8v}=q$ and $q_{\rm 8v}=t$.)
The infinite product-type formulas 
obtained in this paper will be one of the essential ingredients 
in our study of the eight-vertex model in the next paper \cite{S3}.
\bigskip

The plan of the paper is as follows. In Section 2, 
the Macdonald-type difference operator $D$ is introduced, 
and the commutativity $[I(\alpha),D]=0$ is argued.
In Section 3, 
the Weyl group action (denoted by $\pi_m$) on the 
eigenfunctions of $I(\alpha)$ is studied, under the termination condition
for the series $t=q^{m}$ ($m=1,2,3,\cdots$).
In Section 4, it is observed that we have an 
infinite product formula for the eigenfunction of $I(\alpha)$
for the case $(s_1,s_2\cdots,s_n)=(1,t,\cdots,t^{n-1})$.
Comments for the non diagonalizable case 
(which occurs if $s_i=s_j$ is satisfied) are given Section 5. 
A conjecture for the structure of the 
Jordan blocks is obtained.
In Section 6 and Section 7, 
we work with the homogeneous condition
$(s_1,s_2\cdots,s_n)=(1,1,\cdots,1)$.
(Hence $I(\alpha)$ is not diagonalizable.)
The quasi-eigenfunction $F(\alpha)$ is introduced.
Explicit formulas for $F(\alpha)$ are studied for small $n$.
A class of infinite product(-type) formulas for $F(\alpha)$
are conjectured, which takes place at $\alpha=-t^{1/2}$, $\alpha=t$,
$\alpha=\pm q^{1/2}t^{1/2}$, and $\alpha=\pm q^{\ell}t^{1/2}$ ($\ell=1,2,3,\cdots$).
Concluding remarks are given in Section 8.
Appendix A is devoted to the explanation for
the Fock realization of the 
Macdonald difference operators, and difference operators
for the raising operators.
\bigskip

As was in the last paper, 
we use the standard notations
for the $q$-shifted factorials
and the basic hypergeometric series used in 
Gasper and Rahman \cite{GR} (hereafter referred to as GR).
The notation for the elliptic theta function 
$
\Theta_q(z)=(z;q)_\infty(q/z;q)_\infty(q;q)_\infty
$ is used.
\bigskip

\section{Macdonald-type Difference Operator $D$}
In the first paper \cite{S1}, it was conjectured that 
the integral transformation $I(\alpha)$ generates a
commutative family of operators acting on the space ${\cal F}_n$.
In this section, we present another supporting argument for this conjecture.
A Macdonald-type difference operator will be introduced. By examining the
eigenfunctions for $D$, we conjecture that $I(\alpha)$ and $D$ are
commutative.
\bigskip

\subsection{definition of the difference operator $D$}
Let us introduce a difference operator which is
acting on the space of power series ${\cal F}_n$.
\begin{dfn}
Let $s_1,s_2,\cdots,s_n$, $q$ and $t$ be parameters.
Define a difference operator acting on the space 
${\cal F}_n$ by
\begin{eqnarray}
D=D(s_1,\cdots s_n;q,t)
= \sum_{i=1}^n s_i 
\prod_{j<i}\theta_-\left( \zeta_i\over \zeta_j\right)
\prod_{j>i}\theta_+\left( \zeta_j\over \zeta_i\right)
\cdot
T_{q^{-1},\zeta_i}.
\end{eqnarray}
Here $\theta_\pm(\zeta)$ are the series
\begin{eqnarray}
\theta_\pm(\zeta )={1-q^{\pm 1} t^{\mp 1}\zeta \over 1-q^{\pm 1}\zeta}=
1+\sum_{n=1}^\infty (1-t^{\mp 1}) q^{\pm n}\zeta^n,
\end{eqnarray}
and the difference operator $T_{x,\zeta_i}$ is defined by
\begin{eqnarray}
T_{x,\zeta_i}  \cdot g(\zeta_1,\zeta_2,\cdots,\zeta_n)=
g(\zeta_1,\cdots,x\zeta_i,\cdots,\zeta_n).
\end{eqnarray}
\end{dfn}

Note that this difference operator $D$ was derived from the 
Macdonald difference operator $D^1_n$ \cite{Mac},
in the context of the raising operators.
This is explained in Appendix A.
\medskip

\subsection{existence of the eigenfunctions of $D$}
Let us study the existence of the eigenfunctions of the difference operator $D$.
\begin{prop}\label{existing}
Let the parameters $(s_1,s_2,\cdots,s_n)$ and $q$ be generic.
Let $j_1,j_2,\cdots,j_{n-1}$ be nonnegative integers.
In the space ${\cal F}_n$, there exist a unique solution to the equation
\begin{eqnarray}
&&
D f_{j_1,j_2,\cdots,j_{n-1}}(\zeta_1,\cdots,\zeta_n)
=
 \lambda_{j_1,j_2,\cdots,j_{n-1}}
  f_{j_1,j_2,\cdots,j_{n-1}}(\zeta_1,\cdots,\zeta_n),
  \end{eqnarray}
with the conditions
\begin{eqnarray}
&&
f_{j_1,j_2,\cdots,j_{n-1}}(\zeta_1,\cdots,\zeta_n)\\
&=&
\sum_{i_1\geq j_1, \cdots,i_{n-1}\geq j_{n-1}}^\infty
c_{i_1,i_2,\cdots,i_{n-1}}\left(\zeta_2\over \zeta_1\right)^{i_1}
\left(\zeta_3\over \zeta_2\right)^{i_2}
\cdots \left(\zeta_n\over \zeta_{n-1}\right)^{i_{n-1}},\nonumber
\end{eqnarray}
and $c_{j_1,j_2,\cdots,j_{n-1}}=1$,
if and only if 
\begin{eqnarray}
\lambda_{j_1,j_2,\cdots,j_{n-1}}
&=&
\sum_{i=1}^{n}s_i \, q^{-j_{i-1}+j_i} ,\label{eigen-D}
\end{eqnarray}
is satisfied.
Here $j_0=0$ and $j_{n}=0$ are assumed.
\end{prop}

\proof
We have
\begin{eqnarray}
&&D \left(\zeta_2\over \zeta_1\right)^{j_1}\left(\zeta_3\over \zeta_2\right)^{j_2}
\cdots \left(\zeta_n\over \zeta_{n-1}\right)^{j_{n-1}}
=
\left(\zeta_2\over \zeta_1\right)^{j_1}\left(\zeta_3\over \zeta_2\right)^{j_2}
\cdots \left(\zeta_n\over \zeta_{n-1}\right)^{j_{n-1}}\nonumber\\
&&\times
\sum_{i=1}^n s_i \,q^{-j_{i-1}+j_i}
\prod_{j<i}\theta_-\left( \zeta_i\over \zeta_j\right)
\prod_{j>i}\theta_+\left( \zeta_j\over \zeta_i\right), \label{matel-D}
\end{eqnarray}
and it is explicitly seen here that
$D$ is lower triangular
in the monomial basis with respect to the dominance order, or 
the lexicographic order. The diagonal elements are given by Eq.(\ref{eigen-D}). 
If the parameters are generic, 
all the diagonal entries are distinct and
we can construct the eigenfunction.
\qed

By examining the matrix elements of $D$ given by Eq.(\ref{matel-D}),
one finds that
all the eigenfunctions are related by 
shifting the parameters $s_i$.

\begin{prop}\label{eigen-rel}
The eigenfunctions of $D$ 
satisfy
\begin{eqnarray}
&&f_{j_1,j_2,\cdots,j_{n-1}}(\zeta_1,\cdots,\zeta_n)\\
&=&
\prod_{i=1}^{n}
 \zeta_i^{j_{i-1}-j_i}
(T_{q,s_i})^{-j_{i-1}+j_i}\cdot
f_{0,0,\cdots,0}(\zeta_1,\cdots,\zeta_n).\nonumber
\end{eqnarray}
Here, $j_0=0,j_n=0$ are assumed, and $T_{q,s_i}$  denotes
the shift operator acting on the variable $s_i$.
\end{prop}

Note that we have exactly the same property for the
eigenfunctions of the integral transformation $I(\alpha)$.
See Proposition 3.2 of \cite{S1}.

\subsection{commutativity between $I(\alpha)$ and $D$}
Our aim here is to investigate the relationship between the
integral transformation $I(\alpha)$ and the difference operator $D$, 
to claim 
\begin{con}\label{I-vs-D}
The integral transformation $I(\alpha;s_1,\cdots ,s_n,q,t)$ and
the difference operator $D(s_1,\cdots,s_n,q,t)$ are commutative
with each other
\begin{eqnarray}
[I(\alpha;s_1,\cdots,s_n,q,t),D(s_1,\cdots,s_n,q,t)]=0,
\end{eqnarray}
for general $n\geq 2$.
\end{con}

In the first paper \cite{S1}, we obtained the conjecture that 
the integral transformation $I(\alpha)$ generates a commutative family 
of operators acting on the space of series ${\cal F}_n$
(see Conjecture 1.2 of \cite{S1}).
In Appendix A, it is observed that another commutative family of 
Macdonald-type difference operators exists, which are 
acting on ${\cal F}_n$ and 
containing $D(s_1,\cdots,s_n,q,t)$ (see Conjecture \ref{family}).
Therefore, Conjecture \ref{I-vs-D} gives 
us a complementary understanding of the two conjectures for these 
families of commuting operators.
\medskip

In what follow, we argue that the eigenfunctions for 
$I(\alpha)$ and $D$ exactly coincide, 
at least for small $n$ up to certain degree in $\zeta$.
\bigskip

Let us start from the case $n=2$.
\begin{prop}\label{n=2}
The first eigenfunction of $D$
is given by
\begin{eqnarray}
f_{0}(\zeta_1,\zeta_2)=(1-\zeta_2/\zeta_1)\, {}_2\phi_1\left(
{  q t^{-1}, q t^{-1}s_1/s_2 
\atop
q s_1/s_2  }
;q, t\zeta_2/\zeta_1
\right). \label{first-ei}
\end{eqnarray}
All the other eigenfunctions 
$f_{i}(\zeta_1,\zeta_2)$ are given by Proposition \ref{eigen-rel}.
\end{prop}

\proof
Set $f_0=(1-\zeta) g(\zeta)$ and
$g(\zeta)=\sum_{n=0}^\infty g_n \zeta^n$, where $\zeta=\zeta_2/\zeta_1$.
The equation for $f_0$ gives us the difference equation for $g$ as
\begin{eqnarray*}
 s_1 (1-qt^{-1}\zeta)g(q\zeta)+
 s_2 (1-q^{-1}t\zeta)g(q^{-1}\zeta)
 =(s_1+s_2)(1-\zeta)g(\zeta).
\end{eqnarray*}
Solving this with the condition $g_0=1$, we have
\begin{eqnarray*}
g_n={(qt^{-1}, qt^{-1}s_1/s_2;q)_n \over 
(q, q s_1/s_2;q)_n} t^n.
\end{eqnarray*}
\qed

Hence we see
that all the eigenfunctions of $D$ for $n=2$ are
completely identical to the ones for the integral transformation $I(\alpha)$.
(See Theorem 2.1 of \cite{S1}.)
\begin{prop}
The integral transformation $I(\alpha;s_1,s_2,q,t)$ and
the difference operator $D(s_1,s_2,q,t)$ are commutative
on the space ${\cal F}_2$
\begin{eqnarray}
[I(\alpha;s_1,s_2,q,t),D(s_1,s_2,q,t)]=0.
\end{eqnarray}
\end{prop}
\bigskip

Now we proceed to looking at the case $n\geq 3$. By a brute force calculation,
one can observe the following.
\begin{con}\label{n=3}
The first eigenfunction of the difference operator $D$ for the case $n=3$
is given by
\begin{eqnarray}
&&f_{0,0}(\zeta_1,\zeta_2,\zeta_3)\nonumber\\
&=&
\sum_{k=0}^\infty
{
(qt^{-1},qt^{-1},t,t;q)_k \over 
(q,qs_1/s_2,qs_2/s_3,
qs_1/s_3;q)_k}
 (qs_1/s_3)^k (\zeta_3/\zeta_1)^k\label{g-fun}\\
&&\times
\prod_{1\leq i<j\leq 3}
(1-\zeta_j/\zeta_i)
{}_2\phi_1\left(
{  q^{k+1} t^{-1}, q t^{-1}s_i/s_j 
\atop
q^{k+1}s_i/s_j  };q, t\zeta_j/\zeta_i
\right).\nonumber
\end{eqnarray}
\end{con}
We realize the complete coincidence of the
eigenfunctions for $I(\alpha)$ and $D$ for $n=3$ also, at least up to certain 
degrees in $\zeta$.
See Conjecture 3.3 of \cite{S1}. 
\medskip

By using the partial result for $n=4$ given in Proposition 3.4 of \cite{S1},
one can observe that the eigenfunctions for $I(\alpha)$ and $D$
are the same for small degrees in $\zeta$.
Hence Conjecture \ref{I-vs-D} is likely correct for general $n$.
\bigskip

\section{Weyl Group Symmetry}
In this section, we study a hidden symmetry of the eigenfunctions 
of $I(\alpha)$ or $D$ in terms of the Weyl group of type $A_{n-1}$.
This Weyl group symmetry appears when 
the series for the eigenfunctions are terminating, 
at the special points $t=q^{m}$ ($m=1,2,3,\cdots$).
\medskip

Let us fix our notations for the Weyl group action.
Let ${\cal P}_n$ be the space of polynomials in $\zeta_1,\zeta_2,\cdots,\zeta_n$
with coefficients in the field of rational functions in $s_i$'s .
We introduce a representation of the 
Weyl group of type $A_{n-1}$ on the space ${\cal P}_n$ as follows.
\begin{dfn}\label{Weyl}
Let $W(A_{n-1})$ be the Weyl group of type $A_{n-1}$ generated by 
$\sigma_1,\sigma_2,\cdots,\sigma_{n-1}$
with the braid relations $\sigma_i^2={\rm id}$ and
$\sigma_i\sigma_{i+1}\sigma_i=\sigma_{i+1}\sigma_i\sigma_{i+1}$.
Let $m$ be a positive integer.
Define the action of $W(A_{n-1})$ on  ${\cal P}_n$ by
\begin{eqnarray}
&&\pi_m(\sigma_i ) f(\zeta_1,\cdots,\zeta_i,\zeta_{i+1},\cdots,\zeta_n,
s_1,\cdots,s_i,s_{i+1},\cdots,s_n)\\
&=&
\prod_{k=1}^{m-1} {s_{i}-q^{k} s_{i+1} \over s_{i+1}-q^{k} s_i}
f(\zeta_1,\cdots,\zeta_{i+1},\zeta_{i},\cdots,\zeta_n,
s_1,\cdots,s_{i+1},s_{i},\cdots,s_n).\nonumber
\end{eqnarray}
\end{dfn}

Now we are ready to state the main conjecture in this section.
\begin{con}\label{Weyl-con}
If the condition $t=q^{m}$ ($m=1,2,3,\cdots$) is satisfyed,
the first eigenfunction of  $I(\alpha)$ or $D$ for general $n$
becomes terminating
\begin{eqnarray}
\prod_{k=1}^{n}
\zeta_k^{(n-k)m} f_{0,0,\cdots,0}(\zeta_1,\zeta_2,\cdots,\zeta_n)\in{\cal P}_n.
\end{eqnarray}
This polynomial is antisymmetric with respect to the Weyl group action $\pi_m$
\begin{eqnarray}
&&\pi_m(\sigma) \cdot \prod_{k=1}^{n}
\zeta_k^{(n-k)m} f_{0,0,\cdots,0}(\zeta_1,\zeta_2,\cdots,\zeta_n)\nonumber \\
&=&
-\prod_{k=1}^{n}
\zeta_k^{(n-k)m} f_{0,0,\cdots,0}(\zeta_1,\zeta_2,\cdots,\zeta_n)
\qquad ({}^\forall\sigma\in W(A_{n-1})).
\end{eqnarray}
\end{con}
\medskip

In what follows, we present several evidences which support the conjecture.
We will check the termination of the
eigenfunctions and the Weyl group symmetry
under the specialization $t=q^{m}$ ($m=1,2,3,\cdots$), 
by using the explicit formulas for small $n$.
(See Proposition \ref{n=2}, Conjecture \ref{n=3} in this paper, and 
Proposition 3.4 of \cite{S1}.)
\bigskip

Let us consider the case $n=2$.
One finds that the following 
terminating ${}_2\phi_1$ series is $W(A_1)$ symmetric.
\begin{lem}
For $m=1,2,3,\cdots$, the equality
\begin{eqnarray}
&&\pi_m(\sigma_1) \cdot \zeta_1^{m-1}
{}_2\phi_1\left(
{  q t^{-1},  q t^{-1}s_1/s_2
\atop
 q s_1/s_2}
;q, t\zeta_2/\zeta_1
\right)\Biggl|_{t=q^{m}}\nonumber\\
&=&\zeta_1^{m-1}
{}_2\phi_1\left(
{  q t^{-1},  q t^{-1}s_1/s_2
\atop
 q s_1/s_2}
;q, t\zeta_2/\zeta_1
\right)\Biggl|_{t=q^{m}},
\end{eqnarray}
holds for $\sigma_1\in W(A_1)$.
\end{lem}

Hence we have
\begin{prop}
Conjecture \ref{Weyl-con} is true for $n=2$.
\end{prop}
\medskip

Next, let us study the case $n=3$.
One may find a family of $W(A_2)$ symmetric
expressions by using terminating ${}_2\phi_1$ series.
\begin{lem}
Let $m$ be a positive integer and $k$ be a nonnegative integer.
Set
\begin{eqnarray}
\varphi_k &=& 
(\zeta_3/\zeta_1)^k(qs_1/s_3)^k
{(q t^{-1},q t^{-1},t,t;q)_k\over 
(q,qs_1/s_2,qs_2/s_3,qs_1/s_3;q)_k}\\
&&\times
\prod_{1\leq i<j\leq 3}
{}_2\phi_1\left(
{  q^{k+1} t^{-1},  q t^{-1}s_i/s_j
\atop
 q^{{k+1}} s_i/s_j}
;q, t\zeta_j/\zeta_i
\right)\Biggl|_{t=q^{m}}. \nonumber
\end{eqnarray}
Then,
$\zeta_1^{2m-2}\zeta_2^{m-1}\varphi_k$ is terminating and $W(A_2)$ symmetric.
Namely, we have $\zeta_1^{2m-2}\zeta_2^{m-1}\varphi_k \in {\cal P}_3$ and 
\begin{eqnarray}
&&\pi_m(\sigma) \cdot \zeta_1^{2m-2}\zeta_2^{m-1}\varphi_k=
\zeta_1^{2m-2}\zeta_2^{m-1}\varphi_k.
\end{eqnarray}
for ${}^\forall\sigma\in W(A_2)$.
\end{lem}

Since the first eigenfunction for $n=3$ can be written as
a sum $f_{0,0}=\prod_{i<j} (1-\zeta_j/\zeta_i)\sum_k \varphi_k$ (see Conjecture \ref{n=3}), 
we have 
\begin{prop}
Conjecture \ref{Weyl-con} is true for $n=3$, under the assumption that
the formula for  $f_{0,0}(\zeta_1,\zeta_2,\zeta_3)$ given in 
Conjecture \ref{n=3} is correct.
\end{prop}
\medskip

By looking at the partial result for $n=4$ given in Proposition 3.4 of \cite{S1},
one can also observe that the first eigenfunctions for $I(\alpha)$ or $D$
can be decomposed by using 
$W(A_{3})$ antisymmetric terminating series for $t=q^{m}$.
For example, 
we can check the $W(A_{3})$ antisymmetry by
setting $t=q^{m}$ and multiplying
$\zeta_1^{3m}\zeta_2^{2m}\zeta_3^m$ to the combination
\begin{eqnarray*}
&&{\zeta_3\over \zeta_1} \left(q{s_1\over s_3} \right) 
{(q t^{-1})_1(q t^{-1})_1(t)_1(t)_1 \over
(q)_1(qs_{12})_1(qs_{23})_1(qs_{13})_1}
\varphi(1,1,0,1,0,0)\\
&+&{\zeta_4\over \zeta_2} \left(q{s_2\over s_4} \right) 
{(qt^{-1})_1(qt^{-1})_1(t)_1(t)_1 \over
(q)_1(qs_{23})_1(qs_{34})_1(qs_{24})_1}
\varphi(0,1,1,0,1,0)\\
&+&{\zeta_4\over \zeta_1} \left(q{s_1\over s_4} \right) 
{(qt^{-1})_1(qt^{-1})_1(t)_1(t)_1 \over
(q)_1(qs_{12})_1(qs_{24})_1(qs_{14})_1}
\varphi(1,0,0,0,1,1)\nonumber \\
&+& {\zeta_4\over \zeta_1} \left(q{s_1\over s_4} \right) 
{(qt^{-1})_1(qt^{-1})_1(t)_1(t)_1 \over
(q)_1(qs_{13})_1(qs_{34})_1(qs_{14})_1}
\varphi(0,0,1,1,0,1).
\end{eqnarray*}
(As for the notations, see Eq.(66) in \cite{S1}.) 
Thus, we expect that Conjecture \ref{Weyl-con} is correct for $n=4$.

\section{Product Formula for the Eigenfunction}\label{product-f}
In this section,  we study the special case
\begin{eqnarray}
(s_1,s_2,\cdots,s_n)=(1,t,\cdots,t^{n-1}), \label{res-of-s}
\end{eqnarray}
and obtain an infinite product formula for the first
eigenfunction. We
claim the following. 
\begin{con}\label{product_form}
Under the specialization Eq.(\ref{res-of-s}), 
the first eigenfunction of $I(\alpha)$
can be written as the infinite product
\begin{eqnarray}
&&f_{0,0,\cdots,0}(\zeta_1,\zeta_2,\cdots,\zeta_n)
\Biggl|_{(s_1,s_2,\cdots,s_n)=(1,t,\cdots,t^{n-1})} \label{prod-form}\\
&=&
\prod_{1\leq i<j\leq n}(1-\zeta_j/\zeta_i)
 { (qt^{-1}\zeta_j/\zeta_i;q)_\infty \over
 (t\zeta_j/\zeta_i;q)_\infty}.\nonumber
\end{eqnarray}
\end{con}
We give several arguments for this conjecture below.
\bigskip

For $n=2$, we can easily see that 
the infinite series for the eigenfunction $f_0(\zeta_1,\zeta_2)$
reduces into the infinite product. 
\begin{prop}
Conjecture \ref{product_form} is true for $n=2$.
\end{prop}
\proof
Setting $s_1=1,s_2=t$ in Eq.(\ref{first-ei}), and
using the $q$-binomial theorem (Eq. (1.3.2) of GR \cite{GR}), we have
\begin{eqnarray*}
f_0(\zeta_1,\zeta_2)\Bigl|_{(s_1,s_2)=(1,t)}&=&
(1-\zeta)\times
{}_2\phi_1\left(
{  q t^{-1}, qt^{-1}s
\atop
q s }
;q, t\zeta
\right)\Biggl|_{s=t^{-1}}\\
&=&
(1-\zeta) 
{ (qt^{-1}\zeta ;q)_\infty \over
 (t\zeta ;q)_\infty},
\end{eqnarray*}
where $\zeta=\zeta_2/\zeta_1$ and $s=s_1/s_2$.
\qed
\medskip

Next, for $n=3$ we have
\begin{prop}
Conjecture \ref{product_form} is true for $n=3$, under the assumption that
the formula for  $f_{0,0}(\zeta_1,\zeta_2,\zeta_3)$ given in 
Conjecture \ref{n=3} is correct.
\end{prop}
\proof
Changing the order of the summation, and 
using the $q$-Pfaff-Saalsch\"utz formula (Eq. (1.7.2) of GR \cite{GR}), we have
the equality,
\begin{eqnarray*}
&&
\sum_{k=0}^\infty
{(t;q)_k(t;q)_k\over (q;q)_k (qt^{-2};q)_k}(qt^{-2}\zeta)^k
{}_2\phi_1\left({ q^{k+1}t^{-1}, q t^{-3}
\atop q^{k+1}t^{-2}}; q,t\zeta \right)\\
&=&
\sum_{m=0}^\infty 
{(q t^{-1};q)_m(qt^{-3};q)_m \over 
(qt^{-2};q)_m(q;q)_m} t^m \zeta^m 
{}_3\phi_2\left({t,t,q^{-{m}} \atop 
qt^{-1}, q^{-{m}} t^3} ;q,q \right)\\
&=& { (qt^{-1}\zeta;q)_\infty \over
 (t\zeta;q)_\infty}.
\end{eqnarray*}
Then the product formula for the first eigenfunction
\begin{eqnarray}
f_{0,0}(\zeta_1,\zeta_2,\zeta_3)\Biggl|_{(s_1,s_2,s_3)=(1,t,t^2)}=
\prod_{1\leq i<j\leq 3}(1-\zeta_j/\zeta_i)
 { (qt^{-1}\zeta_j/\zeta_i;q)_\infty \over
 (t\zeta_j/\zeta_i;q)_\infty},
 \end{eqnarray}
is derived from Eq.(\ref{g-fun}) by using
the above identity and the $q$-binomial theorem.
\qed
\medskip

For the case $n=4$, we can observe that the degeneration into a product expression
occurs by using the partial result givin in 
Proposition 3.4 of \cite{S1}. 
\medskip

If we assume the commutativity $[I(\alpha),D]=0$, we have
\begin{prop}
Conjecture \ref{product_form} is true for general $n$, under the assumption that
Conjecture \ref{I-vs-D} is correct.
\end{prop}
\proof
It suffices to show
\begin{eqnarray*}
&&D(1,t,\cdots,t^{n-1};q,t)
\prod_{1\leq i<j\leq n}(1-\zeta_j/\zeta_i)
 { (qt^{-1}\zeta_j/\zeta_i;q)_\infty \over
 (t\zeta_j/\zeta_i;q)_\infty}\nonumber \\
 &=&
 (1+t+\cdots+t^{n-1})
 \prod_{1\leq i<j\leq n}(1-\zeta_j/\zeta_i)
 { (qt^{-1}\zeta_j/\zeta_i;q)_\infty \over
 (t\zeta_j/\zeta_i;q)_\infty},
\end{eqnarray*}
for general $n$.
Using the
eigenvalue of the Macdonald difference operator 
$D_n^1$ (see \cite{Mac}), 
we have
\begin{eqnarray*}
&&{\rm LHS}\\
&=&t^{n-1}\prod_{1\leq i<j\leq n}(1-\zeta_j/\zeta_i)
 { (qt^{-1}\zeta_j/\zeta_i;q)_\infty \over
 (t\zeta_j/\zeta_i;q)_\infty}
\sum_{i=1}^n \prod_{j\neq i} {t^{-1} \zeta_i-\zeta_j \over \zeta_i-\zeta_j}
={\rm RHS}.
\end{eqnarray*}
\qed

Thus, it is expected that 
Conjecture \ref{product_form} is true for general $n$.

\section{Generalized Eigenfunctions}
So far, we have been studying the properties of the integral transformation $I(\alpha)$
based on the eigenfunctions. Here in this section, we will make several remarks
for the case $s_i=s_j$. As we will see, the operator $I(\alpha)$ or $D$
acting on ${\cal F}_n$
becomes non diagonalizable, when $n\geq 3$ and
$s_i=s_j$ is satisfied.
Therefore, we need some treatments 
for constructing the generalized eigenfunctions.
In this section, the structure of the eigenspace for 
the operator $I(\alpha)$ or $D$ at
the homogeneous limit
\begin{eqnarray}
s_1=s_2=\cdots=s_n=1,
\end{eqnarray}
will be conjectured. 
\bigskip

By looking at the eigenvalues of $I(\alpha)$ given in Eq.(50) of \cite{S1}
\begin{eqnarray}
\lambda_{j_1,j_2,\cdots,j_{n-1}}(\alpha)
&=&
\prod_{i=1}^{n}
{(\alpha s_i^{-1};q)_{j_{i-1}-j_i} \over 
  (\alpha s_i^{-1} qt^{-1};q)_{j_{i-1}-j_i}},\label{eigenvalues}
\end{eqnarray}
we realize that some of the eigenvalues become degenerate
under the condition: $n\geq 3$ and
$s_i=s_j$ for some $i$ and $j$.
Note, however, that for the case $n=2$, we do not have any degeneration 
of the eigenvalues on the space ${\cal F}_2$.  (So, $I(\alpha)$ remains
diagonalizable, even if $s_1=s_2$.) 
It can be seen from the explicit formulas that
some of the eigenfunctions become divergent,
if we have degenerate eigenvalues. 
(See Proposition 3.2, 
Conjecture 3.3 and Proposition 3.4 in the first paper \cite{S1}.)
In this situation, the integral transformation $I(\alpha)$ (and also for 
the difference operator $D$)
has Jordan blocks.
\medskip

The generalized eigenfunctions can be constructed as follows.
Let $X=X(u)$ be an operator (acting on ${\cal F}_n$) 
depending on a parameter $u$. Assume we have the 
eigenvalues $\lambda(u),\mu(u)$ and the 
eigenfunctions 
$f(u),g(u)$ 
\begin{eqnarray*}
X(u) f(u) =\lambda(u) f(u),\qquad X(u) g(u) =\mu(u) g(u),
\end{eqnarray*}
and the expansions
\begin{eqnarray*}
&&X(u)=X_{[0]}+X_{[1]}u+\cdots,\\
&&\lambda(u)=\lambda_{[0]}+\lambda_{[1]}u+\cdots,\qquad\qquad\quad
\mu(u)=\mu_{[0]}+\mu_{[1]}u+\cdots,\\
&&f(u)={f_{[-1]}\over u} +f_{[0]}+f_{[1]}u+\cdots,\qquad
g(u)=g_{[0]}+g_{[1]}u+\cdots,
\end{eqnarray*}
together with the degeneration condition $\lambda_{[0]}=\mu_{[0]}$.
We also assume that the dimension of the space of
solution to the equation $X_{[0]}f=\lambda_{[0]}f$ is one.
Comparing the coefficients of the above equations in $u$,
we have
\begin{eqnarray*}
&&X_{[0]} f_{[-1]}=\lambda_{[0]}f_{[-1]}, \qquad 
X_{[0]} f_{[0]}+X_{[1]} f_{[-1]}=\lambda_{[0]}f_{[0]}+\lambda_{[1]}f_{[-1]}, \\
&&X_{[0]} g_{[0]}=\mu_{[0]}g_{[0]}, \qquad \quad
X_{[0]} g_{[1]}+X_{[1]} g_{[0]}=\mu_{[0]}g_{[1]}+\mu_{[1]}g_{[0]}, 
\end{eqnarray*}
and so on. Hence we have the proportionality $f_{[-1]}=c g_{[0]}$ and
and the generalized eigenfunctions  as
\begin{eqnarray*}
&& X_{[0]}(f_{[0]}-c g_{[1]})=\lambda_{[0]} (f_{[0]}-c g_{[1]})+
(\lambda_{[1]}-\mu_{[1]}) c g_{[0]}.
\end{eqnarray*}
The generalized eigenfunctions of $I(\alpha)$ (or $D$) can be obtained 
in the above manner.
More degenerate cases such as $s_i=s_j=\cdots =s_k$ can be
treated in a similar manner.
\medskip

Let us introduce some notations.
Let $\Delta=\{\alpha_1,\alpha_2,\cdots,\alpha_{n-1}\}$ be
the set of simple roots for $A_{n-1}$, and $Q_+$ be the 
positive cone of the root lattice. 
We identify the index for the eigenfunctions with
the element in $Q_+$ in the natural way.
We allow to use the same symbol also
for the generalized eigenfunctions, for simplicity. 
We write 
$f_{i,j}(\zeta_1,\zeta_2,\zeta_3)=f_\alpha(\zeta_1,\zeta_2,\zeta_3)$
for $\alpha=i \alpha_1+j \alpha_2$, for example, and
we have the generalized eigenfunction $f_{\alpha_1}$ 
satisfying
\begin{eqnarray*}
I(\alpha)f_{\alpha_1}= 
\lambda_{\alpha_1}f_{\alpha_1}+  \nu_{\alpha_1}(\alpha) f_{\alpha_1+\alpha_2},
\end{eqnarray*}
when $s_2=s_3$ (thus $\lambda_{\alpha_1}=\lambda_{\alpha_1+\alpha_2}$),
and so on.
\medskip

Let us denote the Weyl chamber by 
$C(\Delta)=\{ x|  (x,\alpha_i)\geq 0,\alpha_i\in \Delta \} $.
{}From the explicit formulas for the eigenfunctions
(see Proposition 3.2, Conjecture 3.3 and Proposition 3.4 of \cite{S1}),
we observe the following structures for 
the eigenfunctions of $I(\alpha)$ (and also for $D$)
on the space of series ${\cal F}_n$.
\begin{con}
Let $\alpha \in Q_+$.
If and only if $\alpha \in C(\Delta)$,
the eigenfunction $f_\alpha$ remains finite
at the 
homogeneous limit $s_1=s_2=\cdots=s_n=1$.
Otherwise, $f_{\alpha}$ becomes divergent at the 
homogeneous limit, and a generalized eigenfunction
occurs.
\end{con}
\begin{con}
For the homogeneous case
$s_1=s_2=\cdots=s_n=1$, every eigenspace $V_\lambda$ of 
$I(\alpha)$ or $D$ is uniquely characterized by an element $\alpha \in Q_+$
which is in the Weyl chamber $\alpha \in C(\Delta)$ as
\begin{eqnarray}
V_{\lambda_{\alpha}}=
\bigoplus_{\sigma \in W(A_{n-1}) \atop \sigma(\alpha)\in Q_+} {\bf C} f_{\sigma(\alpha)},
\end{eqnarray}
where $ f_{\alpha}$ denote the (generalized) eigenfunctions.
\end{con}
\begin{con}\label{jordan}
Even in the non diagonalizable cases,
the commutativity among the integral transformations 
$[I(\alpha),I(\beta)]=0$ 
still holds. 
\end{con}
\bigskip

\subsection{example}
To show an example for the explicit formula of the generalized eigenfunctions,
we treat the simplest case.
As we have noted,
the integral transformation $I(\alpha)$ and the difference operator $D$
are diagonalizable on the space of series ${\cal F}_2$, 
even if we spacialise as $s_1=s_2$.
 However,
if we allow negetive powers in $\zeta_2/\zeta_1$ 
and consider the action of $I(\alpha)$ or $D$ on ${\cal F}_2[\zeta_1/\zeta_2]$,
the situation changes and we have to consider
the generalized eigenfunctions at the limit $s_1=s_2$.
\medskip

The eigenfunctions of $I(\alpha)$ or $D$ on the extended space 
${\cal F}_2[\zeta_1/\zeta_2]$
is given by $(\zeta_2/\zeta_1)^i T_{q^{i},s_1}T_{q^{-i},s_2} f_0(\zeta_1,\zeta_2)$ 
($i\in{\bf Z}$),  where $f_0$ is givin in Eq.(\ref{first-ei}).
The explicit formulas for the 
generalized eigenfunctions can be obtained by using the 
above stated method and the
following. 
\begin{lem} 
For $i=0,1,2,\cdots$, 
\begin{eqnarray}
&&{}_2\phi_1\left(
{q t^{-1}, q^{1+2i } t^{-1} s \atop 
q^{1+2i }  s}; q, t \zeta
\right)
=
{}_2\phi_1\left(
{q t^{-1}, q^{1+2i } t^{-1}  \atop 
q^{1+2i}  }; q, t \zeta
\right)\nonumber \\
&&\qquad +(1-s)\sum_{n=0}^\infty 
{(q t^{-1};q)_n (q^{1+2i} t^{-1};q)_n 
\over 
(q;q)_n (q^{1+2i};q)_n } 
(t\zeta )^n\\
&&\qquad \qquad\times
\sum_{k=1}^{n}\left[
{1\over 1-q^{-2i-k }}-{1\over 1-q^{-2i-k}t}\right]+O((1-s)^2),\nonumber
\end{eqnarray}
and for $i=1,2,\cdots$,
\begin{eqnarray}
&&{}_2\phi_1\left(
{qt^{-1}, q^{1-2i } t^{-1} s \atop 
q^{1-2i }  s}; q, t \zeta
\right)\nonumber\\
&=&
{1\over 1-s}
{(q t^{-1};q)_{2i} (q^{1-2i} t^{-1};q)_{2i}
\over 
(q;q)_{2i} (q^{1-2i};q)_{2i-1} } 
(t\zeta)^{2i}
{}_2\phi_1\left(
{q t^{-1}, q^{1+2i } t^{-1}  \atop 
q^{1+2i }  }; q, t \zeta
\right)\nonumber \\
&&\qquad +\sum_{n=0}^{2i-1} 
{(q t^{-1};q)_n (q^{1-2i} t^{-1};q)_n 
\over 
(q;q)_n (q^{1-2i};q)_n } 
(t \zeta)^n\\
&&\qquad+
\sum_{n=2i}^{\infty} 
{(q t^{-1};q)_n (q^{1-2i} t^{-1};q)_n 
\over 
(q;q)_n 
(q^{1-2i};q)_{2i-1} (q;q)_{n-2i} } 
(q \zeta)^n \nonumber\\
&&\qquad\qquad \times \Biggl[ 
\sum_{k=1\atop k\neq 2i}^n {1\over 1-q^{2i-k}}-
\sum_{k=1}^n {1\over 1-q^{2i-k}t}
\Biggr] +O(1-s),\nonumber
\end{eqnarray}
hold.
\end{lem}

\section{Quasi-Eigenfunction of the Integral Transformation $I(\alpha)$}
One of the important characteristics of the basic hypergeometric series 
${}_2\phi_1$ is the existence of 
various infinite product expressions
\begin{eqnarray}
&&{}_2\phi_1\left({a,b\atop b};q,z\right)={(a z;q)_\infty \over (z;q)_\infty}, \label{2phi1-1}\\
&&
{}_2\phi_1\left({a^2,a q\atop a};q,z\right)=(1+a z){(a^2 q z;q)_\infty \over (z;q)_\infty},
\label{2phi1-2}\\
&&
{}_2\phi_1\left({a,-a\atop -q};q,z\right)=
{(a^2 z;q^2)_\infty \over (z;q^2)_\infty},\label{2phi1-3}
\end{eqnarray}
and so on.
We have studied a product formula for the eigenfunctions of $I(\alpha)$
in Section \ref{product-f}. One may find, however, that
possible variety of infinite product formulas for the eigenfunctions 
is not rich enough compared with that of ${}_2\phi_1$.
In this section, we propose a one parameter family of 
`quasi-eigenfunctions' of the integral transformation $I(\alpha)$,
and study it's properties.
It will be pointed out in the next section that a class of 
infinite product(-type) formulas exists for the quasi-eigenfunctions.
\medskip

\subsection{definition of the quasi-eigenfunction $F(\alpha)$}
From now on, we work with the homogeneous condition
\begin{eqnarray}
s_1=s_2=\cdots=s_n=1,
\end{eqnarray}
and we work with the operator
$I(\alpha)=I(\alpha;1,1,\cdots,1,q,t)$.
\medskip

We introduce the
`quasi-eigenfunction  
$F(\alpha)$' as follows.
\begin{dfn} \label{Falpha}
Set the parameters as $s_1=s_2=\cdots=s_n=1$.
The quasi-eigenfunction  
$F(\alpha)=F(\zeta_1,\zeta_2,\cdots,\zeta_n;\alpha,q,t)$
is defined by the `covariant transformation property'
\begin{eqnarray}
{\rm (I)}&&
I(\alpha q^{-1}t ) \cdot F(\alpha)=F(\alpha q^{-1}t).
\end{eqnarray}
and the `initial condition'
\begin{eqnarray}
{\rm (II)}&&
F(t^{1/2})
=\prod_{1\leq i<j\leq n}(1-\zeta_j/\zeta_i)
{(qt^{-1/2}\zeta_j/\zeta_i;q)_\infty \over 
(t^{1/2}\zeta_j/\zeta_i;q)_\infty }.
\end{eqnarray}
\end{dfn}

Here, an explanation is in order.
First, $F(\alpha)$ can be constructed 
at $\alpha=t^{1/2},q^{-1}t^{3/2},q^{-2}t^{5/2},\cdots$ 
by 
the iterative action of $I(\alpha)$ as
\begin{eqnarray}
&&F(q^{-1}t^{3/2})=I(q^{-1}t^{3/2})\cdot F(t^{1/2}), \label{it-ac-1}\\
&&F(q^{-2}t^{5/2})=I(q^{-2}t^{5/2})\cdot I(q^{-1}t^{3/2})\cdot F(t^{1/2}),\label{it-ac-2}
\end{eqnarray}
and so on. Then the function $F(\alpha)$ is obtained by the analytic continuation
with respect to the parameter $\alpha$
from these discrete points. 
One can check that the 
above definition of $F(\alpha)$ is well defined 
for small degrees in $\zeta$, 
by performing explicit analytic continuation of the coefficients.
\bigskip

It is expected that 
the function $F(\alpha)$
satisfies another transformation property 
which is similar to the condition $\rm (I)$.
\begin{con}\label{equiv-tr}
The function $F(\alpha)$ satisfies the 
transformation property
\begin{eqnarray}
%
{\rm (I')}&&
I(\alpha^{-1} q)\cdot  F(\alpha)=F(\alpha q^{-{1}}t).
\end{eqnarray}
\end{con}
This can be checked up to certain order.
Note taht, the both transformation properties 
$\rm (I)$ and $\rm (I')$ will be needed in the next paper, while we study the 
vertex operator of the eight-vertex model \cite{S3}.

\subsection{explicit formula of $F(\alpha)$ for $n=2$}
For the case $n=2$, one can easily obtain an explicit formula of  $F(\alpha)$.
\begin{prop}\label{Falpha-n=2}
Let $F(\alpha)$ be the following ${}_2\phi_1$ (or ${}_4\phi_3$) series:
\begin{eqnarray}
&&F(\alpha)=F(\zeta_1,\zeta_2;\alpha,q,t)\nonumber\\
&=&
(1-\zeta_2/\zeta_1){}_2\phi_1\left(
{ q t^{-1},\alpha q t^{-1}\atop
\alpha^{-1}  q}; q,\alpha^{-1} t \zeta_2/\zeta_1
 \right)\label{Falpha-2}\\
&=&{}_4\phi_3\left(
{q t^{-{1\over 2}},-q t^{-{1\over 2}},
t^{-1}, \alpha  t^{-1}\atop
t^{-{1\over 2}},-t^{-{1\over 2}},
\alpha^{-1} q}; q,\alpha^{-1} t \zeta_2/\zeta_1
 \right)\nonumber.
\end{eqnarray}
Then this $F(\alpha)$ satisfies 
the transformation property $\rm (I)$ and 
the initial condition $\rm (II)$ in Definition \ref{Falpha} for $n=2$.
\end{prop}

\proof
The initial condition $\rm (II)$ can be checked
by using the $q$-binomial theorem (Eq. (1.3.2) of GR \cite{GR}):
\begin{eqnarray*}
F(t^{1/2})
=(1-\zeta_2/\zeta_1)
{(q t^{-1/2}\zeta_2/\zeta_1;q)_\infty \over 
(t^{1/2}\zeta_2/\zeta_1;q)_\infty }.
\end{eqnarray*}
Next, let $f_i(\zeta_1,\zeta_2)$'s be the eigenfunctions of $I(\alpha)$
given in Eq.(17) of \cite{S1}, and
\begin{eqnarray}
\lambda_i(\alpha)
&=&
{(\alpha^{-1} s_1 t ;q)_i \over (\alpha^{-1} s_1q;q)_i}
{(\alpha s_2^{-1};q)_i \over 
 (\alpha s_2^{-1}q t^{-1};q)_i}
\left(q t^{-1} \right)^i, 
\end{eqnarray}
be the corresponding eigenvalues (see Eq.(16) of \cite{S1}).
One finds that the series Eq.(\ref{Falpha-2}) 
can be neatly expanded in terms of the
eigenfunctions $f_i(\zeta_1,\zeta_2)$ 
(see  Lemma \ref{lem--2} given below).
Therefore the transformation property $\rm (I)$
can be shown as
\begin{eqnarray*}
&& I(\alpha q^{-1}t )\cdot F(\zeta_1,\zeta_2;\alpha,q,t)\\
&=&
\sum_{i=0}^\infty f_i(\zeta_1,\zeta_2) \lambda_i(\alpha q^{-1}t )
{(\alpha;q)_i \over (\alpha^{-1}  q;q)_i}
{
( q t^{-1}, q^{i+1}t^{-1};q)_i\over 
(q , q^{i} ;q)_i}
\alpha^{-i}t^i \\
&=&
F(\zeta_1,\zeta_2;\alpha q^{-1}t,q,t).
\end{eqnarray*}
Note that the other transformation property $\rm (I')$ given in
Conjecture \ref{equiv-tr} can also be checked.
\qed

It remains to examine the 
the eigenfunction expansion of $F(\alpha)$.
This is accomplished by using a summation formula for ${}_6\phi_5$ series.
\begin{lem}\label{lem--2}
Set $s_1=s_2=1$.
With respect to the eigenfunctions $f_i(\zeta_1,\zeta_2)$ of $I(\alpha)$, 
the function $F(\zeta_1,\zeta_2;\alpha,q,t)$ 
given in Eq.(\ref{Falpha-2}) is expanded as
\begin{eqnarray}
&&F(\zeta_1,\zeta_2;\alpha,q,t)\\
&=&
\sum_{i=0}^\infty f_i(\zeta_1,\zeta_2)
{(\alpha;q)_i \over (\alpha^{-1}  q;q)_i}
{
(q t^{-1}, q^{i+1}t^{-1};q)_i\over 
(q , q^{i} ;q)_i}
\alpha^{-i}t^i .\nonumber
\end{eqnarray}
\end{lem}

\proof
From the explicit formula (given in Eq.(17) of \cite{S1}) for $s_1=s_2=1$
\begin{eqnarray*}
f_j(\zeta_1,\zeta_2)=
\zeta^j\times
{}_4\phi_3\left(
{   q^{2j}t^{-1} ,  q^{j+1} t^{-{1\over 2}},
-  q^{j+1} t^{-{1\over 2}} ,t^{-1}
\atop
  q^{j}t^{-{1\over 2}} ,-  q^{j}t^{-{1\over 2}},
  q^{2j+1}}
;q, t\zeta
\right),
\end{eqnarray*}
we have the expansion with respect to the monomials
\begin{eqnarray*}
f_j(\zeta_1,\zeta_2)&=&\sum_{i=0}^\infty c_{ij}(\zeta_2/\zeta_1)^i, \\
c_{ij}&=&
\left\{
\begin{array}{ll}
\displaystyle
{
(q^{2j}t^{-1}, q^{j+1}t^{-{1\over 2}},
 -q^{j+1}t^{-{1\over 2}},t^{-1};q)_{i-j}
\over 
(q^{2j+1},q^{j}t^{-{1\over 2}},
 -q^{j}t^{-{1\over 2}},q;q)_{i-j}}
t^{i-j} & (i\geq j) ,\\[4mm]
0 & (i<j) ,
\end{array}
\right.
\end{eqnarray*}
where we have 
denoted $\zeta=\zeta_2/\zeta_1$ (see Eq.(30) of \cite{S1}). 
Writing 
\begin{eqnarray*}
F(\zeta_1,\zeta_2;\alpha,q,t)=
\sum_{i=0}^\infty a_i (\zeta_2/\zeta_1)^i=
\sum_{i=0}^\infty f_i(\zeta_1,\zeta_2) b_i,
\end{eqnarray*}
we have the matrix equation
$a =C b$, where $a={^t}(a_0,a_1,\cdots)$, $b={^t}(b_0,b_1,\cdots)$  and
$C=(c_{ij})_{0\leq i,j<\infty}$ 
is the lower triangular matrix defined as above. Therefore
we have $b=C^{-1}a$.

In Proposition 2.7 of \cite{S1}, we have obtained 
the entries for the inverse $C^{-1}=(d_{ij})_{0\leq i,j,\infty}$ as
\begin{eqnarray}
d_{ij}&=&
\left\{
\begin{array}{ll}
\displaystyle
{
(q^{i+j+1}t^{-1},t;q)_{i-j}
\over 
(q^{i+j},q;q)_{i-j}} & (i\geq j) ,\\[4mm]
0 & (i<j) .
\end{array}
\right.
\end{eqnarray}
{}From the expression in Eq.(\ref{Falpha-2}),
the coefficients $a_i$'s read 
\begin{eqnarray}
a_j=
{(qt^{-{1\over 2}},-qt^{-{1\over 2}},
t^{-1},\alpha t^{-1};q)_j
\over 
(t^{-{1\over 2}},- t^{-{1\over 2}},
\alpha^{-1}  q,q;q)_j} (\alpha^{-1}t)^j.\nonumber
\end{eqnarray}
Hence we arrive at the formula for $b_i$'s  as follows
\begin{eqnarray*}
&&b_i=\sum_{j=0}^id_{ij}a_j\\
&=&
{(t, q^{i+1}t^{-1};q)_{i}
\over
(q, q^{i};q)_{i}}
\sum_{j=0}^i 
{(q^{-i}, q^{i};q)_{j}
\over
(q^{-i+1}t^{-1}, q^{i+1}t^{-1};q)_{j}}
(qt^{-1})^j
a_j\\
&=&
{(t,q^{i+1}t^{-1};q)_{i}
\over
(q,q^{i};q)_{i}}\\
&&\times
{}_6\phi_5
\left(
{
qt^{-{1\over 2}},-qt^{-{1\over 2}},
q^{-i},q^{i},t^{-1},\alpha t^{-1}
\atop
t^{-{1\over 2}},-t^{-{1\over 2}},
q^{i+1}t^{-1},q^{-i+1}t^{-1},\alpha^{-1} q}
;q,\alpha^{-1}q
\right)\\
&=&
{(t,q^{i+1}t^{-1};q)_{i}
\over
(q,q^{i};q)_{i}}
{(qt^{-1},\alpha^{-1} q^{-i+1};q)_{i}
\over
(q^{-i+1}t^{-1},\alpha^{-1}q;q)_{i}}\\
&=&
{(\alpha;q)_i \over (\alpha^{-1} q;q)_i}
{
(qt^{-1},q^{i+1}t^{-1};q)_i\over 
(q,q^{i} ;q)_i}
(\alpha^{-1}t)^i .
\end{eqnarray*}
Here we have used the summation formula for
the ${}_6\phi_5$ (Eq. (2.4.2) of GR \cite{GR}). 
\qed
\bigskip

\subsection{conjectural form of $F(\alpha)$ for $n=3$}
Working with the monomial basis,
one can perform
the iterative actions of $I(\alpha)$ (at lease for small degrees in $\zeta$)
as in Eq.(\ref{it-ac-1}), (\ref{it-ac-2}) etc., and
study the analytic continuation with respect to the variable $\alpha$.
Then an explicit expression of the quasi-eigenfunction $F(\alpha)$
for $n=3$ is guessed as follows.
\begin{con}\label{Falpha-3}
The quasi-eigenfunction for $n=3$ is written as
\begin{eqnarray}
&&F(\zeta_1,\zeta_2,\zeta_3;\alpha,q,t)\label{F(alphan)=3}\\
&=&
\sum_{k=0}^\infty
{
(\alpha^{-2}t,qt^{-1},qt^{-1};q)_k \over 
(q,\alpha^{-1}q,\alpha^{-1}q;q)_k}
(q\zeta_3/\zeta_1)^k
{}_2\phi_1 \left({\alpha^{-1}, q^{-{k}} \atop \alpha q^{-k+1}};
q, \alpha t\right)
\nonumber\\
&&\times
\prod_{1\leq i<j\leq 3}
(1-\zeta_j/\zeta_i)\;
{}_2\phi_1\left(
{  q^{k+1} t^{-1}, \alpha q t^{-1}
\atop
\alpha^{-1}q^{{k+1}} };q, \alpha^{-1}t\zeta_j/\zeta_i
\right).\nonumber
\end{eqnarray}
\end{con}

When $F(\alpha)$ is expanded in terms of the generalized 
eigenfunctions, we observe the following.
\begin{con}\label{embed}
We have the embedding
\begin{eqnarray}
F(\alpha)\in \bigoplus_{k=0}^\infty  V_{\lambda_{k \theta}},
\end{eqnarray}
where $\theta=\alpha_1+\alpha_2$ denotes the maximal root for $A_2$.
\end{con}
This means that 
the initial condition $\rm (II)$ given in Definition \ref{Falpha}
is quite restrictively chosen.

\subsection{partial result for $n=4$}
For $n=4$, a brute force calculation gives us
the following observation.
\begin{con}\label{Falpha-n=4}
The quas-eigenfunction 
$F(\zeta_1,\zeta_2,\zeta_3,\zeta_4;\alpha,q,t)$
is given by the series
\begin{eqnarray}
&&F(\zeta_1,\zeta_2,\zeta_3,\zeta_4;\alpha,q,t)\\
&=&
\sum_{k=0}^\infty Y_{k,k,0}+\sum_{k=0}^\infty Y_{k,k+1,1}
+\sum_{k=0}^\infty Y_{k+1,k+1,1}+\cdots, \nonumber
\end{eqnarray}
on the subspace of ${\cal F}_4$ spanned by the monomials
\begin{eqnarray*}
\left(\zeta_2\over \zeta_1\right)^{i_1}
\left(\zeta_3\over \zeta_2\right)^{i_2}
\left(\zeta_4\over \zeta_3\right)^{i_3}\qquad 
(0 \leq i_1<\infty,0 \leq i_2<\infty,0 \leq i_3\leq 1).
\end{eqnarray*}
Here the series $Y_{i,j,k}$ are defined as follows:
\begin{eqnarray}
Y_{k,k,0}
&=&q^{k}\left(\zeta_3\over \zeta_1\right)^k
{
(\alpha^{-2}t)_k(qt^{-1})_k(qt^{-1})_k \over 
(q)_k(\alpha^{-1}q)_k(\alpha^{-1}q)_k}
{}_2\phi_1 \left({\alpha^{-1}, q^{-{k}} \atop \alpha q^{-k+1}};
q, \alpha t\right)\nonumber\\
&&\times
\phi(k,k,0,k,0,0),\qquad\qquad ({\rm for}\,\,k=0,1,2,\cdots),\\
Y_{0,1,1}
&=&q\left(\zeta_4\over \zeta_2\right)
{
(\alpha^{-2}t)_1(qt^{-1})_1(qt^{-1})_1 \over 
(q)_1(\alpha^{-1}q)_1(\alpha^{-1}q)_1}
{}_2\phi_1 \left({\alpha^{-1}, q^{-{1}} \atop \alpha};
q, \alpha t\right)\nonumber\\
&&\times
\phi(0,1,1,0,1,0),\\
Y_{1,1,1}
&=&q\left(\zeta_4\over \zeta_1\right)
{
(\alpha^{-2}t)_1(qt^{-1})_1(qt^{-1})_1 \over 
(q)_1(\alpha^{-1}q)_1(\alpha^{-1}q)_1}
{}_2\phi_1 \left({\alpha^{-1}, q^{-{1}} \atop \alpha};
q, \alpha t\right)\nonumber\\
&&\times
[\phi(1,0,0,0,1,1)+\phi(0,0,1,1,0,1)],\nonumber\\
&-&q
\left(\zeta_4\over \zeta_1\right)
{
(\alpha^{-2}t)_1(qt^{-1})_1(qt^{-1})_1
(qt^{-1})_1 \over 
(q)_1
(\alpha^{-1}q)_1(\alpha^{-1}q)_1
(\alpha^{-1}q)_1}\\
&&\times
{}_2\phi_1 \left({\alpha^{-1}, q^{-{1}} \atop \alpha };
q, \alpha t\right)
\left[
{(1-\alpha)q\over (1-q)\alpha}\,
{}_2\phi_1 \left({\alpha^{-1}, q^{-{1}} \atop \alpha };
q, \alpha t\right) \right.\nonumber\\
&&+ \left.
{\alpha-q\over (1-q)\alpha}\,
{}_2\phi_1 \left({\alpha^{-1}, q^{-{1}} \atop \alpha};
q, \alpha qt\right)
\right]\nonumber\\
&&\times
\phi(1,1,1,1,1,1),\nonumber
\end{eqnarray}
\begin{eqnarray}
Y_{k,k+1,1}
&=&-q\left(\zeta_3\over \zeta_1\right)^k
\left(\zeta_4\over \zeta_2\right)
{
(\alpha^{-2}t)_k(qt^{-1})_k(qt^{-1})_k
(qt^{-1})_1 \over 
(q)_{k-1}(q)_1
(\alpha^{-1}q)_k(\alpha^{-1}q)_k
(\alpha^{-1}q)_1}\nonumber\\
&&\times
{}_2\phi_1 \left({\alpha^{-1}, q^{-{1}} \atop \alpha };
q, \alpha t\right)
\left[
{(1-\alpha)q^{k}\over (1-q^{k})\alpha}\,
{}_2\phi_1 \left({\alpha^{-1}, q^{-{k}} \atop \alpha q^{-k+1}};
q, \alpha t\right) \right.\nonumber\\
&&+ \left.
{\alpha-q^{k}\over (1-q^{k})\alpha}\,
{}_2\phi_1 \left({\alpha^{-1}, q^{-{k}} \atop \alpha q^{-k+1}};
q, \alpha qt\right)
\right]\nonumber\\
&&\times
\phi(k,k,1,k,1,1)\nonumber\\
&+&q\left(\zeta_3\over \zeta_1\right)^k
\left(\zeta_4\over \zeta_2\right)
{
(\alpha^{-2}t)_{k+1}(qt^{-1})_{k+1}(qt^{-1})_k
(qt^{-1})_1 \over 
(q)_{k}(q)_1
(\alpha^{-1}q)_{k+1}(\alpha^{-1}q)_k
(\alpha^{-1}q)_1}\\
&&\times
{}_2\phi_1 \left({\alpha^{-1}, q^{-{1}} \atop \alpha };
q, \alpha t\right)
{}_2\phi_1 \left({\alpha^{-1}, q^{-{k}} \atop \alpha q^{-k+1}};
q, \alpha qt\right)\nonumber\\
&&\times
\phi(k,k+1,1,k,1,1),\qquad ({\rm for}\,\,k=1,2,3,\cdots),\nonumber
\end{eqnarray}
and
\begin{eqnarray}
&&Y_{k+1,k+1,1} \nonumber\\
&=&
-q\left(\zeta_3\over \zeta_1\right)^k
\left(\zeta_4\over \zeta_1\right)
{
(\alpha^{-2}t)_{k+1}(qt^{-1})_{k+1}(qt^{-1})_{k+1}
(qt^{-1})_1 \over 
(q)_{k}(q)_1
(\alpha^{-1}q)_{k+1}(\alpha^{-1}q)_{k+1}
(\alpha^{-1}q)_1}\nonumber\\
&&\times
{}_2\phi_1 \left({\alpha^{-1}, q^{-{1}} \atop \alpha };
q, \alpha t\right)
\left[
{(1-\alpha)q^{k+1}\over (1-q^{k+1})\alpha}\,
{}_2\phi_1 \left({\alpha^{-1}, q^{-{k+1}} \atop \alpha q^{-k}};
q, \alpha t\right)\right.\nonumber\\
&&+\left.
{\alpha-q^{k+1}\over (1-q^{k+1})\alpha}\,
{}_2\phi_1 \left({\alpha^{-1}, q^{-{k+1}} \atop \alpha q^{-k}};
q, \alpha qt\right)
\right]\nonumber\\
&&\times
\phi(k+1,k+1,1,k+1,1,1)\nonumber\\
&-&
q\left(\zeta_3\over \zeta_1\right)^k
\left(\zeta_4\over \zeta_1\right)
{
(\alpha^{-2}t)_{k}(qt^{-1})_{k}(qt^{-1})_k
(qt^{-1})_1 \over 
(q)_{k-1}(q)_1
(\alpha^{-1}q)_{k}(\alpha^{-1}q)_k
(\alpha^{-1}q)_1}\\
&&\times
{}_2\phi_1 \left({\alpha^{-1}, q^{-{1}} \atop \alpha };
q, \alpha t\right)
\left[
{(1-\alpha)q^{k}\over (1-q^{k})\alpha}\,
{}_2\phi_1 \left({\alpha^{-1}, q^{-{k}} \atop \alpha q^{-k+1}};
q, \alpha t\right)\right.\nonumber\\
&&+\left.
{\alpha-q^{k}\over (1-q^{k})\alpha}\,
{}_2\phi_1 \left({\alpha^{-1}, q^{-{k}} \atop \alpha q^{-k+1}};
q, \alpha qt\right)
\right]\nonumber\\
&&\times
\phi(k,k,1,k,1,1)\nonumber\\
&+&q\left(\zeta_3\over \zeta_1\right)^k
\left(\zeta_4\over \zeta_1\right)
{
(\alpha^{-2}t)_{k+1}(qt^{-1})_{k+1}(qt^{-1})_k
(qt^{-1})_1 \over 
(q)_{k}(q)_1
(\alpha^{-1}q)_{k+1}(\alpha^{-1}q)_k
(\alpha^{-1}q)_1}\nonumber\\
&&\times
{}_2\phi_1 \left({\alpha^{-1}, q^{-{1}} \atop \alpha };
q, \alpha t\right)
{}_2\phi_1 \left({\alpha^{-1}, q^{-{k}} \atop \alpha q^{-k+1}};
q, \alpha qt\right)\nonumber\\
&&\times
[\phi(k,k,1,k+1,1,1)+\phi(k+1,k,1,k,1,1)],\qquad ({\rm for}\,\,k=1,2,3,\cdots),\nonumber
\end{eqnarray}
were, we have used the notations 
$(a)_k=(a;q)_k$ and
\begin{eqnarray}
&&\phi(k_{12},k_{23},k_{34},k_{13},k_{24},k_{14})\\
&=&
\prod_{1\leq i<j\leq 4}
(1-\zeta_j/\zeta_i)
\cdot {}_2\phi_1\left(
{  q^{k_{ij}+1} t^{-1},\alpha q t^{-1} 
\atop
\alpha^{-1}q^{{k_{ij}+1}}  };q, \alpha^{-1}t\zeta_j/\zeta_i
\right).\nonumber
\end{eqnarray} 
 \end{con}

\section{Product Formulas for the Quasi-Eigenfunction $F(\alpha)$}
The aim of this section is 
to present several conjectures for
infinite product(-type) formulas for the quasi-eigenfunction $F(\alpha)$
at some particular values of $\alpha$.
One realizes that these can be regarded 
as multi variable generalizations of the product formulas for ${}_2\phi_1$ given in
Eqs.(\ref{2phi1-1})-(\ref{2phi1-3}).

\subsection{case $\alpha=- t^{1/2}$ and $\alpha=t$}
We observe that our quasi-eigenfunction $F(\alpha)$ 
neatly factorizes into product expressions if we specialize
the parameter $\alpha$ to $-t^{1/2}$ or $t$.
\begin{con}\label{degenF}
The infinite product formulas
\begin{eqnarray}
&&
F(- t^{1/2})
=\prod_{1\leq i<j\leq n}(1-\zeta_j/\zeta_i)
{(- qt^{-1/2}\zeta_j/\zeta_i;q)_\infty \over 
(- t^{1/2}\zeta_j/\zeta_i;q)_\infty }, \label{degen1}\\
&&
F( t )
=\prod_{1\leq i<j\leq n\atop {\rm step} 2}(1-\zeta_j/\zeta_i)
{( qt^{-1 }\zeta_j/\zeta_i;q)_\infty \over 
( t \zeta_j/\zeta_i;q)_\infty }, \label{degen2}
\end{eqnarray}
hold.
Here, we have used the notation
 $\prod_{1\leq i<j\leq n\atop {\rm step} 2} f_{ij}=
f_{13}f_{15}\cdots f_{24}f_{26}\cdots$.
\end{con}

For the case $n=2$, we can prove the factorization.
\begin{prop}
Conjecture \ref{degenF} is true for $n=2$.
\end{prop}
\proof
{}From the explicit formula Eq.(\ref{Falpha-2}) and the
$q$-binomial theorem, we have
\begin{eqnarray*}
&&
F(- t^{1/2})
=(1-\zeta_2/\zeta_1)
{(- qt^{-1/2}\zeta_2/\zeta_1;q)_\infty \over 
(- t^{1/2}\zeta_2/\zeta_1;q)_\infty }, \qquad
F( t )
=1.
\end{eqnarray*}
\qed
\medskip

Next, we prove the factorization for the conjectural 
expression of $F(\alpha)$ for $n=3$.
\begin{prop}
Conjecture \ref{degenF} is true for $n=3$, under the assumption that
the formula for  $F(\alpha)$ given in 
Conjecture \ref{Falpha-3} is correct.
\end{prop}
\proof
Setting $\alpha=-t^{1/2}$ in Eq.(\ref{F(alphan)=3}),
the summand becomes zero for $k>0$. Then using
the $q$-binomial theorem, we have
\begin{eqnarray*}
&&
F(- t^{1/2})
=\prod_{1\leq i<j\leq 3}(1-\zeta_j/\zeta_i)
{(- qt^{-1/2}\zeta_j/\zeta_i;q)_\infty \over 
(- t^{1/2}\zeta_j/\zeta_i;q)_\infty }.
\end{eqnarray*}
 
We proceed to proving the case $\alpha=t$. 
First, note that we have 
\begin{eqnarray*}
{}_2\phi_1\left(
{  q^{k+1} t^{-1}, \alpha q t^{-1}
\atop
\alpha^{-1}q^{{k+1}} };q, \alpha^{-1}t\zeta
\right)\Biggl|_{\alpha=t}={1\over 1-\zeta}.
\end{eqnarray*}
Next, using Jackson's transformation formula
(Eq.(iii) in Exercise 1.15 of GR \cite{GR})
\begin{eqnarray*}
{}_2\phi_1\left({q^{-n},b\atop c};q,z\right)
={(c/b;q)_n\over (c;q)_n} 
{}_3\phi_2\left({q^{-n},b,bzq^{-n}/c\atop bq^{1-n}/c,0};q,q\right),
\end{eqnarray*}
we have
\begin{eqnarray*}
{}_2\phi_1 \left({\alpha^{-1}, q^{-{k}} \atop \alpha q^{-k+1}};
q, \alpha t\right)\Biggl|_{\alpha=t}
=q^{-k} t^k {(t^{-2};q)_k\over (t^{-1};q)_k}{1+t^{-1}q^{k}\over 1+t^{-1}}.
\end{eqnarray*}
Therefore from Eq.(\ref{F(alphan)=3}) we have 
\begin{eqnarray*}
F(\alpha)\Bigl|_{\alpha=t}
&=&
\sum_{k=0}^\infty 
{(t^{-1};q)_k\over (q;q)_k}
 {(t^{-2};q)_k\over (t^{-1};q)_k}{1+t^{-1}q^{k}\over 1+t^{-1}}(t \zeta_3/\zeta_1)^k\\
&=&
{1\over 1+t^{-1}}\left(
{(t^{-1}\zeta_3/\zeta_1;q)_\infty \over (t\zeta_3/\zeta_1;q)_\infty}+
t^{-1}
{(qt^{-1}\zeta_3/\zeta_1;q)_\infty \over (qt\zeta_3/\zeta_1;q)_\infty}
\right)\\
&=&(1-\zeta_3/\zeta_1){(qt^{-1}\zeta_3/\zeta_1;q)_\infty
\over (t\zeta_3/\zeta_1;q)_\infty}.
\end{eqnarray*}
\qed

For the case $n=4$, one can check the factorization by using the partial result
given in Conjecture \ref{Falpha-n=4}.
We expect that Conjecture \ref{degenF} is true for general $n$.
\medskip

\subsection{case $\alpha=\pm q^{1/2} t^{1/2}$}
In the case $\alpha=\pm q^{1/2} t^{1/2}$, we observe that $F(\alpha)$ 
can be written as a multiple of a Pfaffian and an infinite product.
\begin{con}\label{prod-pfaff}
Let $n$ be a positive even integer.
At $\alpha=q^{1/2} t^{1/2}$, 
the expression 
\begin{eqnarray}
&&
F(q^{1/2} t^{1/2}) \nonumber  \\
&=&
{\rm Pfaffian}
\left( {1-\zeta_j^2/\zeta_i^2 \over 
(1-q^{-1/2}t^{1/2}\zeta_j/\zeta_i)(1-q^{1/2}t^{-1/2}\zeta_j/\zeta_i)}
 \right)_{1\leq i,j\leq n}\\
 &&\times
 \prod_{1\leq i<j\leq n}
{(q^{1/2}t^{-1/2}\zeta_j/\zeta_i;q)_\infty \over 
(q^{1/2}t^{1/2}\zeta_j/\zeta_i;q)_\infty }, \nonumber
\end{eqnarray}
hold. The formula for odd $n$ is obtained by taking the limit
$\zeta_n\rightarrow 0$.
\end{con}
Note the case $\alpha=-q^{1/2} t^{1/2}$ is
obtained by negating $t^{1/2}$ in the above formula.
\medskip

Let us examine the above statement for small $n$.
\begin{prop}
Conjecture \ref{prod-pfaff} is true for $n=2$.
\end{prop}
\proof
{}From the explicit formula Eq.(\ref{Falpha-2}) and the
product formula Eq.(\ref{2phi1-2}), we have
\begin{eqnarray}
F(q^{1/2} t^{1/2}) 
&=&
 {1-\zeta_2^2/\zeta_1^2 \over 
(1-q^{-1/2}t^{1/2}\zeta_2/\zeta_1)(1-q^{1/2}t^{-1/2}\zeta_2/\zeta_1)}\\
 &&\times
 {(q^{1/2}t^{-1/2}\zeta_2/\zeta_1;q)_\infty \over 
(q^{1/2}t^{1/2}\zeta_2/\zeta_1;q)_\infty }. \nonumber
\end{eqnarray}
\qed

Next, we prove the product-Pfaffian formula for the conjectural 
expression of $F(\alpha)$ for $n=3$.
\begin{prop}
Conjecture \ref{prod-pfaff} is true for $n=3$, under the assumption that
the formula for  $F(\alpha)$ given in 
Conjecture \ref{Falpha-3} is correct.

\end{prop}
\proof
Setting $\alpha=q^{1/2}t^{1/2}$ in Eq.(\ref{F(alphan)=3}),
the summand becomes zero for $k\geq 2$. Note that we have 
\begin{eqnarray*}
\!\!\!\!&&{}_2\phi_1\left(
{  q^{k+1} t^{-1}, \alpha q t^{-1}
\atop
\alpha^{-1}q^{{k+1}} };q, \alpha^{-1}t\zeta
\right)\Biggl|_{\alpha=q^{1/2}t^{1/2} \atop k=0{~~~~~~~}}=
(1+\zeta){(q^{3/2}t^{-1/2}\zeta;q)_\infty \over 
(q^{-1/2}t^{-1/2}\zeta;q)_\infty },\\
\!\!\!\!&&{}_2\phi_1\left(
{  q^{k+1} t^{-1}, \alpha q t^{-1}
\atop
\alpha^{-1}q^{{k+1}} };q, \alpha^{-1}t\zeta
\right)\Biggl|_{\alpha=q^{1/2}t^{1/2}\atop k=1{~~~~~~~}}=
{(q^{3/2}t^{-1/2}\zeta;q)_\infty \over 
(q^{-1/2}t^{-1/2}\zeta;q)_\infty },
\end{eqnarray*}
from
the $q$-binomial theorem and the product formula Eq.(\ref{2phi1-2}). 
Using these and simplifying the
rational factor in front of the infinite product, we have
\begin{eqnarray}
F(q^{1/2} t^{1/2}) 
&=&
\left( {1-\zeta_2^2/\zeta_1^2 \over 
(1-q^{-1/2}t^{1/2}\zeta_2/\zeta_1)(1-q^{1/2}t^{-1/2}\zeta_2/\zeta_1)}\right.  \nonumber\\
&& +{1-\zeta_3^2/\zeta_2^2 \over 
(1-q^{-1/2}t^{1/2}\zeta_3/\zeta_2)(1-q^{1/2}t^{-1/2}\zeta_3/\zeta_2)}\\
&&- \left.
 {1-\zeta_3^2/\zeta_1^2 \over 
(1-q^{-1/2}t^{1/2}\zeta_3/\zeta_1)(1-q^{1/2}t^{-1/2}\zeta_3/\zeta_1)}
\right) \nonumber\\
 &&\times
 \prod_{1\leq i<j\leq 3}
{(q^{1/2}t^{-1/2}\zeta_j/\zeta_i;q)_\infty \over 
(q^{1/2}t^{1/2}\zeta_j/\zeta_i;q)_\infty }. \nonumber
\end{eqnarray}
\qed 

The case  $n=4$ can be examined by using the partial result
given in Conjecture \ref{Falpha-n=4}.

\subsection{case $\alpha=\pm q^{\ell} t^{1/2}$}
Let us move on to the case 
$\alpha=\pm q^{\ell} t^{1/2}$ ($\ell=1,2,\cdots$).
For $\alpha=\pm q^{1/2} t^{1/2}$, it was argued that the Pfaffian
appears in front of the infinite product. For 
$\alpha=\pm q^{\ell} t^{1/2}$ ($\ell=1,2,\cdots$),
we have another rational expression.
\medskip 

Introduce some notations to describe the structure of the rational factor.
Let $J$ be the index set
$
J=\{\sigma=(\sigma_1,\sigma_2,\cdots,\sigma_{\ell})|\sigma_i=\pm\}
$, 
and introduce $2^{\ell}\times 2^{\ell}$ matrices
\begin{eqnarray}
G_\ell(\zeta )=
\left( \gamma_{\ell,\sigma,\sigma'}(\zeta)
\right)_{\sigma,\sigma'\in J}, \label{GspinB}
\end{eqnarray}
by the following recursive rule.
\begin{dfn}
Set the entries of the matrix $G_1(z)$ as
\begin{eqnarray}
&&\gamma_{1,+,+}(\zeta)=1,\\
&&\gamma_{1,+,-}(\zeta)=
{1-q t^{-1}\zeta \over 1-\zeta}{1+t^{1/2}\zeta \over 1+qt^{-1/2} \zeta},\\
&&\gamma_{1,-,+}(\zeta)=
{1-q^{-1}t \zeta \over 1-\zeta}{1+t^{-1/2}\zeta \over 1+q^{-1}t^{1/2} \zeta},\\
&&\gamma_{1,-,-}(\zeta)=1.
\end{eqnarray}
Then define the matrices $G_\ell(z)$ for $\ell\geq2$ 
recursively by
\begin{eqnarray}
G_{\ell}(\zeta)&=&
\left( 
\begin{array}{cc}
G_{\ell-1} (\zeta) & \gamma_{1,+,-}(\zeta)G_{\ell-1}(\zeta q)\\
\gamma_{1,-,+}(\zeta)G_{\ell-1}(\zeta q^{-1})&G_{\ell-1} (\zeta)
\end{array}
\right).\label{GspinBend}
\end{eqnarray}
\end{dfn}
Using this notation, our observation can be stated as follows.
\begin{con}\label{degenF-2}
Let $n\geq 2$ denote the number of variables for 
the quasi-eigenfunction $F(\alpha)$.
Let $\ell$ be a positive integre, and $G_\ell(\zeta)$ be the matrix of 
rational functions in $\zeta$ given 
as above.
Set 
\begin{eqnarray}
&&\mu_\sigma=
\prod_{i=1}^{\ell}\left(q^{(i-1)/2}t^{1/4}\right)^{\sigma_i},
\end{eqnarray}
for $\sigma=(\sigma_1,\sigma_2,\cdots,\sigma_\ell)\in J$.
Then the following expression for the
quasi-eigenfunction for $\alpha=-q^{\ell}t^{1/2}$ holds:
\begin{eqnarray}
&&F(-q^{\ell}t^{1/2})\nonumber\\
&=&
\left(\sum_{\sigma \in J}\mu_\sigma\right)^{-n}
\prod_{1\leq i<j\leq n}(1-\zeta_j/\zeta_i)
{(- qt^{-1/2}\zeta_j/\zeta_i;q)_\infty \over 
(- t^{1/2}\zeta_j/\zeta_i;q)_\infty }\\
&&\times
\sum_{\sigma_1,\sigma_2,\cdots,\sigma_n\in J}
\prod_{i=1}^n \mu_{\sigma_i}
\cdot
\prod_{1\leq i<j\leq n} \gamma_{\ell,\sigma_i,\sigma_j}(\zeta_j/\zeta_i).\nonumber
\end{eqnarray}
\end{con}
Note the case $\alpha=q^{\ell} t^{1/2}$ is
obtained by negating $t^{1/2}$ in the above formula.
\medskip

By a tedious explicit calculation, one can observe the 
validity of the above conjecture for small $n$.

\subsection{case $\alpha=-1$}
Finally, let us make a short comment on the case $\alpha=-1$.
By looking at the explicit formula of $F(\alpha)$ given in Eq.(\ref{Falpha-2}),
we have 
\begin{prop}For $n=2$, we have the infinite product expression
\begin{eqnarray}
&&
F(-1)
=(1-\zeta_2/\zeta_1)
{(-q^2t^{-1}\zeta_2/\zeta_1;q^2)_\infty \over 
(-t\zeta_2/\zeta_1;q^2)_\infty }.
\end{eqnarray}
\end{prop}
\proof This follows from the product formula Eq.(\ref{2phi1-3}). 
\qed
\medskip

A simple product formula may `not' easily be found for 
$\alpha=-1$ if $n\geq 3$. 
Nevertheless, we are interested in the series $F(-1)$ for general $n$,
since it is related with the study of the
eight-vertex model in a nontrivial manner.
We will argue this in some detail in the next paper \cite{S3}.

\section{Concluding Remarks}
In this paper,
we have studied several properties for the integral transformation $I(\alpha)$,
which was introduced in the first paper \cite{S1}.
\medskip

It was observed that $I(\alpha)$ and the Macdonald-type difference operator
$D$ are commutative with each other (Conjecture \ref{I-vs-D}).
The commutativity $[I(\alpha),D]=0$ was rather 
unexpected, at least to the present author. 
We hope that 
the relationship between the eight-vertex model and the 
Macdonald polynomials will be investigated further.
\medskip

Using the explicit formulas for the eigenfunctions, it was observed that 
the series for the eigenfunction becomes terminating 
under the conditions $t=q^{m}$ ($m=1,2,\cdots$).
For such cases, we found 
the Weyl group symmetry as a hidden symmetry (Conjecture \ref{Weyl-con}).
This Weyl group symmetry helps us when we 
try to construct hypergeometric-type formula
for the eigenfunction for general $q$ and $t$.
\medskip

We have introduced the quasi-eigenfunction $F(\alpha)$ (Definition \ref{Falpha}),
using the action of  $I(\alpha)$ for the homogeneous 
limit $s_1=\cdots=s_n=1$.
The series $F(\alpha)$ has some resemblance to
the first eigenfunction $f_{0}$.
First, if we disregard the dependence on the 
`dynamical' parameter $\alpha$ 
in the transformation property 
$I(\alpha q^{-1}t ) \cdot F(\alpha)=F(\alpha q^{-{1}}t)$, one 
can compare this with the equation
$I(\alpha)\cdot f_{0}=\lambda_{0}f_{0}$.
Next, we have similar hypergeometric-type series expressions
for $f_0$ and $F(\alpha)$ (for $n=3$, see Conjecture \ref{n=3} and
Conjecture \ref{Falpha-3}).
Note also that the infinite product formula for the
first eigenfunction (Conjecture \ref{product_form}) and
the initial condition $\rm (II)$ for $F(\alpha)$ given in 
Definition \ref{Falpha} look similar to each other.
\medskip

For $n=2$, we have studied $F(\alpha)$ in terms of the
eigenfunctions $f_i$ (Proposition \ref{Falpha-n=2}).
For $n\geq 3$ (and $s_1=\cdots=s_n=1$), 
however, $I(\alpha)$ is non diagonalizable 
and our study based on the generalized eigenfunctions
becomes very much complicated.
Nevertheless, one may find
some nontrivial observation as Conjecture \ref{embed}.
\medskip

In Section 7, we have obtained a variety of infinite product-type formulas
for the quasi-eigenfunction $F(\alpha)$.
These were checked for small $n$ by using
the explicit formulas and several transformation and summation
formulas for the basic hypergeometric series.
\medskip

In the next paper \cite{S3}, 
the matrix elements of the
vertex operators for the eight-vertex model will be studied in some detail,
based on the product(-type) formulas for $F(\alpha)$ obtained in this paper.
Our aim there is
to construct a class of Heisenberg representations which
gives us the following description for the vertex operators:
\begin{eqnarray}
\langle \Phi(\zeta_1)\Phi(\zeta_2)\cdots \Phi(\zeta_n)\rangle 
=
\prod_{1\leq i<j\leq n}
{\xi(\zeta_j^2/\zeta_i^2;p,q) \over 1-\zeta_j/\zeta_i} \cdot
F(-1,p^{1/2},q),
\end{eqnarray} 
where
\begin{eqnarray}
&&
\xi(z;p,q)={(q^2z;p,q^4)_\infty (pq^2z;p,q^4)_\infty \over 
(q^4z;p,q^4)_\infty (pz;p,q^4)_\infty }.
\end{eqnarray}Note that
the basic parameters here have been switched to the ones for
the eight-vertex model as  $q\rightarrow p^{1/2}=p_{\rm 8v}^{1/2}$ and 
$t\rightarrow q=q_{\rm 8v}$.
\medskip

Let us give some examples.
Set $p^{1/2}=q^{3/2}$. From the 
transformation property $\rm (I)$ in Definition \ref{Falpha}, and 
the product formula given in Conjecture \ref{degenF},
we have the $n$-fold integral representation for $F(-1,q^{3/2},q)$ as
\begin{eqnarray}
F(-1,q^{3/2},q)=I(-1,q^{3/2},q) \cdot F(-q^{1/2},q^{3/2},q).
\end{eqnarray}
In the same way, we have the integral representations 
$F(-1,-q^2,q)=I(-1,-q^2,q) \cdot F(q,-q^2,q)$ 
for $p^{1/2}=-q^2$, and 
$F(-1,q^3,q)=I(-1,q^3,q) \cdot F(-p^{1/4}q^{1/2},q^3,q)$ for $p^{1/2}=q^3$
(see product formulas in Conjecture \ref{degenF} and 
Conjecture \ref{prod-pfaff}).
These three cases have already been discussed in the previous paper \cite{S}.
\medskip

If we set $p^{1/2}=q^{(2\ell+1)/2\ell}$ ($\ell=1,2,\cdots$), 
we obtain the $\ell \times n$-fold integral representation of
$F(-1)=F(-1,q^{(2\ell+1)/2\ell},q)$ by 
$I(\alpha)=I(\alpha,q^{(2\ell+1)/2\ell},q)$ as
\begin{eqnarray}
F(-1)=I(-1)\cdots I(-pq^{-3/2})\cdot I(-p^{1/2}q^{-1/2}) \cdot F(-q^{1/2}).
\end{eqnarray}
In this way, 
one may construct various integral representations for $F(-1)$ 
by using the infinite product(-type) formulas given in
Conjecture \ref{degenF},
Conjecture \ref{prod-pfaff} and Conjecture \ref{degenF-2}.
In the next paper, we will construct a class of Heisenberg representations of 
the vertex operator $\Phi(\zeta)$, which are consistent with
the conjectures for $F(\alpha)$ obtained in this article.
\bigskip

\noindent
{\it Acknowledgment.}~~~
This work is supported by the  Grant-in-Aid for Scientific Research 
(C) 16540183.
The author thank Yasushi Kajihara for pointing out
the work by Lassalle and Schlosser \cite{LS}.

\appendix
\section{Macdonald Difference Operators}
In this appendix, we revisit 
the Heisenberg representation of the Macdonald difference operators
$D_n^r$
and the Macdonald symmetric function $Q_\lambda(x;q,t)$.
In this appendix, we use the standard notations for the
Macdonald polynomials \cite{Mac}. Note that $q$ and $t$
in this appendix are the ordinary parameters for the 
Macdonald polynomials, and should not be confused with the ones for 
the eight-vertex model $p_{\rm 8v}$ and $q_{\rm 8v}$.

\subsection{commuting operators on the Fock space}
The Macdonald difference operators (see \cite{Mac}) acting on the
ring of symmetric polynomials $\Lambda_{n,F}$ ($F={\bf Q}(q,t)$) are defined by
\begin{eqnarray}
&&D^r_n=\sum_{I\subset \{1,2,\cdots,n\}\atop |I|=r} A_I(x;t)\prod_{i\in I}T_{q,x_i},
\label{D1}\\
&&A_I(x;t)=t^{r(r-1)/2}\prod_{i\in I\atop j\not\in I} {t x_i-x_j \over x_i-x_j}.
\end{eqnarray}
These operators are commutative with each other $[D_n^r,D_n^s]=0$.
Our aim in this appendix is to treat
the commutative family generated by the Macdonald difference operators
over the space of symmetric functions 
$\Lambda_{F}$, namely, in the infinitely many variable situation.
One finds that a use of the Heisenberg algebra will make 
our argument simple.
\bigskip

Introduce the Heisenberg algebra generated by $a_n$ ($n\in{\bf Z}_{\neq 0}$)
satisfying the commutation relations
\begin{eqnarray}
[a_m,a_n]=m{1-q^{|m|} \over 1-t^{|m|}}\delta_{m+n,0}.
\end{eqnarray}
Let $|0\rangle $ be the vacuum vector
satisfying
$a_n |0\rangle=0$ ($n=1,2,\cdots$),
and ${\cal F}$ be the Fock space
\begin{eqnarray}
{\cal F}=F[a_{-1},a_{-2},\cdots]|0\rangle.
\end{eqnarray}
We have the natural identification between the Fock space ${\cal F}$ and the
ring of the symmetric functions $\Lambda_F$ by the rule:
\begin{eqnarray}
{\cal F}\simeq \Lambda_F:\quad
a_{-n} \longleftrightarrow p_n,\qquad |0\rangle\longleftrightarrow 1.
\end{eqnarray}
Here $p_n$ denotes the power sum function $p_n=\sum_i x_i^n$.
\medskip

In  \cite{AMOS} (see also \cite{S-Lec}), 
it was shown that the Macdonald difference operator $E$ defined by
\begin{eqnarray}
&&E = \lim_{\longleftarrow}E_n:\Lambda_F\longrightarrow \Lambda_F,\label{E1}\\
&&E_n = t^{-n}D_n^1 -\sum_{i=1}^n t^{-i}
:\Lambda_{n,F}\longrightarrow \Lambda_{n,F},
\end{eqnarray}
can be realized in terms of the Heisenberg algebra.
It reads,
\begin{eqnarray}
&&H_1\equiv(t-1)E+1=\eta_0,\label{eta0}\\
&&\eta(z)=\sum_{n\in{\bf Z}} \eta_n z^{-n}=\;\; : 
\exp \left( -\sum_{n\neq 0} {1-t^n\over n}a_n z^{-n}\right):,
\end{eqnarray}
where the symbol $:\cdots:$ means the usual normal ordering of the oscillators.
\medskip

The above construction can be extended, and
one can obtain a commutative family of operators $H_r$ acting on the 
Fock space ${\cal F}$.
Let us introduce the $H_r$'s as follows.
\begin{dfn}
Define the operators $H_r$ for $r=1,2,3,\cdots$ by
\begin{eqnarray}
H_r={[r]_{t^{-1}}!\over n!}
\oint_{C_1} {dz_1\over 2 \pi i z_1}\cdots \oint_{C_r}  {dz_r\over 2 \pi i z_r} 
\prod_{1\leq i<j\leq r} \omega(z_j/z_i)
: \eta(z_1)\eta(z_2)\cdots\eta(z_r):, \label{E_r}
 \end{eqnarray}
where the contours $C_i$ are circles 
$|z_i|=1$, and 
$\omega(z)$ is defined by
\begin{eqnarray}
\omega(z)=
{(1-z)(1-z^{-1}) \over (1-t^{-1}z)(1-t^{-1}z^{-1})}=
{2\over 1+t^{-1}}+\sum_{k=1}^\infty {t^{-k}-t^{-k+1}\over 1+t^{-1}}(z^k+z^{-k}).
\end{eqnarray}
\end{dfn}
Note that we have normalized $H_r$ to simplify our later discussion,
using the notations
\begin{eqnarray}
[n]_x={1-x^{n}\over 1-x},\qquad [n]_x!=\prod_{k=1}^n[k]_x.
\end{eqnarray}

Then we claim the following.
\begin{prop}\label{H-H}
On the bosonic Fock space ${\cal F}$, the commutation relations
\begin{eqnarray}
[H_r,H_s]=0,
\end{eqnarray}
hold, for all $r$ and $s$.
\end{prop}
In the next section, we give a proof based on the commutativity of the
Macdonald difference operators $D_n^r$.

\subsection{proof of $[H_r,H_s]=0$.}
We prove Proposition \ref{H-H}, by using several 
operator product formulas and the commutation reltation $[D_n^r,D_n^s]=0$.
\bigskip

Let us introduce the generating function for the
Macdonald symmetric function of length one,
\begin{eqnarray}
\phi(x)
=\sum_{n=0}^\infty \phi_{-n} x^n
=\;\; : \exp \left( \sum_{n= 1}^\infty {1\over n} {1-t^n\over 1-q^n}a_{-n} x^{n}\right):.
\end{eqnarray}
Namely we have $\phi_{-n} \bigl|_{a_{-n}\rightarrow p_n}=Q_{(n)}(q,t)$.
The following formulas will be used 
for studying the action of $H_r$.
\begin{prop}\label{phi-prop}
We have the operator product formula
\begin{eqnarray}
&&\eta(z)\phi(x)=\mu(x/z):\eta(z)\phi(x):, \label{OPE-eta-phi}\\
&&\mu(z)={1-z\over 1-t z}=1+\sum_{k=1}^\infty (t^{k}-t^{k-1})z^k, 
\end{eqnarray}
the difference property
\begin{eqnarray}
:\eta(tx)\phi(x):|0\rangle&=&
 \phi(qx)|0\rangle, \label{dif-prop}
\end{eqnarray}
and the expansion in terms of symmetric polynomials in $x_i$'s
\begin{eqnarray}
&&\phi(x_1)\phi(x_2)\cdots \phi(x_n)|0\rangle 
=
\sum_{\lambda,\ell(\lambda)\leq n} 
\phi_{-\lambda} m_\lambda(x)|0\rangle . \label{expan}
\end{eqnarray}
In Eq.(\ref{expan}), the summation is taken over the partitions, the notation 
$\phi_{-\lambda}=\phi_{-\lambda_1}\phi_{-\lambda_2}\cdots \phi_{-\lambda_n}$
is used for $\lambda=(\lambda_1,\cdots,\lambda_n)$,
and the monomial symmetric polynomial in $x_1,\cdots,x_n$ is denoted by
$m_\lambda(x)$.
\end{prop}

By using Proposition \ref{phi-prop}, 
one can study the action of the operators $H_r$
on the vector $\phi(x_1)\cdots \phi(x_n)|0\rangle$.
\begin{prop}\label{H_r}
The equality
\begin{eqnarray}
&&H_r \cdot  \phi(x_1)\cdots \phi(x_n)|0\rangle \label{eqHr}\\
&=&
t^{-r n}\sum_{k=0}^r (t-1)^k \left[r\atop k\right]_t [k]_t! D_n^k \cdot
\phi(x_1)\cdots \phi(x_n)|0\rangle , \nonumber
\end{eqnarray}
holds,
where $\left[n\atop m \right]_t =[n]_t!/[m]_t![n-m]_t!$, and
$D_n^k$ is the Macdonald difference operator
acting on the variable $x_i$'s
given in Eq.(\ref{D1}).
\end{prop}

In the usual notation, we have $\phi_{-\lambda}|0\rangle=g_\lambda(x;q,t)$
(under the identification ${\cal F}\simeq \Lambda_F$).
Since the functions $g_\lambda(x;q,t)$ form 
the dual basis to the monomial basis $(m_\lambda)$,
the vectors $\phi_{-\lambda}$ form 
a basis of the Fock space ${\cal F}$. Therefore, from
the commutation relations $[D_n^r,D_n^s]=0$,
Proposition \ref{H_r} means that the operators $H_r$ 
are mutually commutative on the Fock space. This proves Proposition \ref{H-H}.
It remains to show Proposition \ref{H_r}.

\proof
First, we note the identity
\begin{eqnarray}
{\rm Symm} \prod_{1\leq i<j\leq n} {1-x_j/x_i \over 1-t^{-1}x_j/x_i}=
{[n]_{t^{-1}}!\over n!} \prod_{1\leq i<j\leq n} \omega(x_j/x_i), \label{ide}
\end{eqnarray}
where the symbol `Symm' means the symmetrization 
in the variavle $x_i$'s.
By using this, we have the expression
\begin{eqnarray*}
&&H_r \cdot \phi(x_1)\cdots \phi(x_n)|0\rangle\\
&=&
\oint_{C_1} \cdots \oint_{C_r}
{d w_1\over 2\pi i w_1}\cdots {d w_r\over 2\pi i w_r}
\prod_{i=1}^r \prod_{j=1}^n {w_i-x_j\over w_i-t x_j}\\
&&\times
\prod_{1\leq i<j\leq n} {w_i-w_j \over w_i-t^{-1}w_j}
:\eta(w_1)\cdots \eta(w_r)\phi(x_1)\cdots\phi(x_n) :|0\rangle,
\end{eqnarray*}
were the integration contours $C_i$ enclose poles at 
$w_i=t x_1,t x_2,\cdots,t x_n$,
$w_i=t^{-1}w_{i+1},t^{-1}w_{i+2},\cdots,t^{-1}w_{r}$ and $w_i=0$.
This multiple integral can be easily performed and 
written as a multiple summation.
Then we symmetrize the summation by
using the identity Eq.(\ref{ide}) again in the form
\begin{eqnarray*}
{\rm Symm} \prod_{1\leq i<j\leq m} {t x_i -x_j \over x_i-x_j}=
{[m]_{t}!\over m!}.
\end{eqnarray*} 
After the symmetrization,
we can organize the result in terms of the
Macdonald difference operators $D_n^r$, and obtain Eq.(\ref{eqHr}). \qed

\subsection{difference equation for the raising operator}
Using the Heisenberg representation of the Macdonald difference operators,
we obtain a difference equation for the
raising operator of the Macdonald Polynomial.
\bigskip

Let $s_1,s_2,\cdots,s_n$ be parameters, and
consider the space of formal power series
$F[[x_2/x_1,x_3/x_2,\cdots,x_n/x_{n-1}]]$.
Here and hereafter, we will work with $F={\bf Q}(q,t,s_1,s_2,\cdots,s_n)$.
We introduce modified Macdonald difference operators as follows.
\begin{dfn}
Define the difference operators 
$D^r(s_1,s_2,\cdots,s_n,q,t)$
acting on $F[[x_2/x_1,x_3/x_2,\cdots,x_n/x_{n-1}]]$ by
\begin{eqnarray}
&&D^r(s_1,s_2,\cdots,s_n,q,t)\\
&=&
\sum_{I\subset \{1,2,\cdots,n\}\atop |I|=r} 
\prod_{i\in I\atop {j\not\in I \atop j<i}}\theta_-\left(x_i\over x_j\right)
\prod_{i\in I\atop {j\not\in I \atop j>i}}\theta_+\left(x_j\over x_i\right)
\prod_{i\in I}s_iT_{q^{-1},x_i},\nonumber
\end{eqnarray}
where $\theta_\pm(x)$ are the series
\begin{eqnarray}
\theta_\pm(x)={1-q^{\pm 1} t^{\mp 1}x\over 1-q^{\pm 1}x}=
1+\sum_{n=1}^\infty (1-t^{\mp 1}) q^{\pm n}x^n. \label{thetapm}
\end{eqnarray}
\end{dfn}

In Section 2, we have studied $D^1(s_1,s_2,\cdots,s_n,q,t)$ 
(note we have 
denoted $D=D^1$ for simplicity), and 
proved the existence of the eigenfunctions.
Let us take the eigenfunction 
\begin{eqnarray}
&&D^1(s_1,s_2,\cdots,s_n,q,t) \cdot f=
\sum_{i=1}^n s_i\, f, \label{Dif}
\end{eqnarray}
with the normalization $f=1+\cdots$. 
Then our claim for
the integral representation or the raising operator for the Macdonald 
symmetric function is stated as follows.
\begin{prop}
Let $f$ be as above.  
Let $\lambda=(\lambda_1,\lambda_2,\cdots,\lambda_n)$
be a partition and 
$Q_{\lambda}(q,t)$
the Macdonald symmetric function (given as a vector in the Fock space).
Then we have
\begin{eqnarray}
&&Q_{\lambda}(q,t)=
\oint \cdots \oint {dx_1\over 2\pi i x_1}\cdots {dx_n\over 2\pi i x_n}
x_1^{-\lambda_1}\cdots x_n^{-\lambda_n}\label{Integral} \\
&&\qquad\qquad \times
f(x_1,\cdots,x_n) 
\phi(x_1)\cdots \phi(x_n)|0\rangle, \nonumber
\end{eqnarray}
by specializing the parameters in $f$ as $s_i=t^{n-i}q^{\lambda_i}$. Here the
integrals $\oint {dx_k\over 2\pi i x_k}$ mean to take the constant term in $x_k$.
\end{prop}

\proof
By using Proposition \ref{H_r} we have
\begin{eqnarray*}
&&H_1 \cdot(\mbox{RHS of Eq.(\ref{Integral})})\\
&=&
 \oint \cdots \oint {dx_1\over 2\pi i x_1}\cdots {dx_n\over 2\pi i x_n}
x_1^{-\lambda_1}\cdots x_n^{-\lambda_n}f(x_1,\cdots,x_n) \\
&&\times (t^{-n}(t-1) D_n^1+t^{-n})
\phi(x_1)\cdots \phi(x_n)|0\rangle\\
&=&
 \oint \cdots \oint {dx_1\over 2\pi i x_1}\cdots {dx_n\over 2\pi i x_n}
x_1^{-\lambda_1}\cdots x_n^{-\lambda_n}
\phi(x_1)\cdots \phi(x_n)|0\rangle \\
&&
\times \Biggl[
t^{-n}(t-1) D^1(t^{n-1}q^{\lambda_1},t^{n-2}q^{\lambda_1},
\cdots,q^{\lambda_n},q,t) +t^{-n}\Biggr]\\
&&\times f(x_1,\cdots,x_n)\\
&=&\left((t-1)\sum_{i=0}^\infty t^{-i} (q^{\lambda_i}-1)+1\right)
(\mbox{RHS of Eq.(\ref{Integral})}).
\end{eqnarray*} 
Here we have made suitable $q^{-1}$-shifts  in $x_i$ to transform the action of 
the operator 
$D_n^1$ on the vector $\phi(x_1)\cdots \phi(x_n)|0\rangle$ to 
the one on the series $f(x_1,\cdots,x_n)$. Note that one has to reinterpret the
rational factors in $D_n^1$ as series acting on the space of symmetric polynomials
$\Lambda_{n,{\bf Q}(q,t)}$ by using the series Eq.(\ref{thetapm}). 
Then we arrive at the difference operator 
$D^1(t^{n-1}q^{\lambda_1},t^{n-2}q^{\lambda_1},
\cdots,q^{\lambda_n},q,t)$ which acts on the series $f(x_1,\cdots,x_n)$.
Since $f(x_1,\cdots,x_n)$ satisfies the difference equation Eq.(\ref{Dif}),
we obtain ${\rm RHS}\propto Q_{\lambda}(q,t)$. The proportionality
can be argued
by using the Pieri formula \cite{Mac}.
\qed

For example,
if $n=2$, the difference equation Eq.(\ref{Dif}) gives us
\begin{eqnarray}
f(x_1,x_2)=(1-x_2/x_1){}_2\phi_1 
\left({ qt^{-1}, qt^{-1}s_1/s_2 \atop q s_1/s_2};q, t x_2/x_1\right).
\end{eqnarray}
By setting $s_1=t q^{\lambda_1},s_2=q^{\lambda_2}$,
we recover Jing and J\'ozefiak's result \cite{JJ}.
For the case $n=3$, the solution to the difference equation Eq.(\ref{Dif}) 
can be guessed as (see Conjecture \ref{n=3})
\begin{eqnarray}
&&f(x_1,x_2,x_3)\nonumber\\
&=&
\sum_{k=0}^\infty
{
(qt^{-1},qt^{-1},t,t;q)_k \over 
(q,q s_1/s_2,q s_2/s_3,
q s_1/s_3;q)_k}
 (q s_1/s_3)^k (x_3/x_1)^k\label{n=3-Mac}\\
&&\times
\prod_{1\leq i<j\leq 3}
(1-x_j/x_i)\;
{}_2\phi_1\left(
{  q^{k+1} t^{-1},q t^{-1}s_i/s_j ,
\atop
q^{k+1}s_i/s_j  };q, t x_j/x_i
\right).\nonumber
\end{eqnarray}
We set $s_1=t^2 q^{\lambda_1},s_2=t q^{\lambda_2},s_3=q^{\lambda_3}$,
for the calculation of $Q_{(\lambda_1,\lambda_2,\lambda_3)}$.
\medskip

Recently, Lassalle and Schlosser found a general formula 
for the raising operator of the Macdonald polynomials \cite{LS}.
They found a way to invert the Pieri formula, and derived 
the raising operator in terms of certain determinant expressions.
It was checked that the series in 
Eq.(\ref{n=3-Mac}) for $n=3$ 
agrees with their determinant formula up to certain degree in $x_i$'s.
\bigskip

It is not a straight forward task to derive 
the commutativity of the operator
$D^r(s_1,s_2,\cdots,s_n,q,t)$'s
from  $[H_r,H_s]=0$, since the expression Eq.(\ref{Integral})
has a kernel and does not specify $f$ uniquely.
Nevertheless, it is naturally expected that 
the $D^r(s_1,s_2,\cdots,s_n,q,t)$'s are commutative with each other.
\begin{con}\label{family}
On the space $F[[x_2/x_1,x_3/x_2,\cdots,x_n/x_{n-1}]]$, we have
\begin{eqnarray}
[D^r(s_1,s_2,\cdots,s_n,q,t),D^s(s_1,s_2,\cdots,s_n,q,t)]=0,
\end{eqnarray}
for all $r$ and $s$.
\end{con}

\end{document}